\documentclass{article}
\usepackage{amssymb,amsmath}
\usepackage{graphics}
   
\def\qed{\hfill {\hbox{${\vcenter{\vbox{               
   \hrule height 0.4pt\hbox{\vrule width 0.4pt height 6pt
   \kern5pt\vrule width 0.4pt}\hrule height 0.4pt}}}$}}}
  
\def\tr{\triangleright}

\def\bar{\overline}  

\newenvironment{proof}[1][Proof]{\smallskip\noindent{\bf #1.}\quad}%
{\qed\par\medskip}

\newtheorem{theorem}{Theorem}
\newtheorem{definition}{Definition}

\newtheorem{proposition}[theorem]{Proposition}
\newtheorem{corollary}[theorem]{Corollary}
\newtheorem{conjecture}{Conjecture}
\newtheorem{example}{Example}
\newtheorem{remark}{Remark}

\author{{\begin{tabular}{c} Sam Nelson \\
\small{\texttt{knots@esotericka.org}}\end{tabular}}
\and {\begin{tabular}{c} John Vo\\
\small{\texttt{johnvo06@yahoo.com}}\end{tabular}}
\and
\small{\begin{tabular}{c}
Department of Mathematics, University of California, Riverside\\
900 University Avenue, Riverside, CA, 92521\end{tabular} }}

\date{}

\title{\Large \textbf{Matrices and Finite Biquandles}}

\begin{document}

\maketitle

\begin{abstract}
We describe a way of representing finite biquandles with $n$ elements as 
$2n\times 2n$ block matrices. Any finite biquandle defines an invariant
of virtual knots through counting homomorphisms. The counting
invariants of non-quandle biquandles can reveal information not present in
the knot quandle, such as the non-triviality of the virtual trefoil and 
various Kishino knots. We also exhibit an oriented virtual knot which is 
distinguished from both its obverse and its reverse by a finite biquandle 
counting invariant. We classify biquandles of order 2, 3 and 4 and provide 
a URL for our Maple programs for computing with finite biquandles.
\end{abstract}

\textsc{Keywords:} Finite biquandles, virtual knots, knot invariants

\textsc{2000 MSC:} 57M25, 57M27

\section{\large\textbf{Introduction}}

Much recent work has been done on the knot invariants defined by counting
homomorphisms from a knot quandle into a finite target quandle (\cite{I},
\cite{DL}, etc.). A \textit{quandle} is the algebraic structure
obtained by assigning a generator to each arc in an oriented knot diagram 
with a binary operation at each crossing -- specifically, set 
$c=a\tr b$ when $b$ is the overcrossing arc, $a$ the arc on the right-hand 
side of $b$ when looking in the positive direction of the oriented\footnote{
This use of orientation is not strictly necessary; the definition in \cite{J}
uses unoriented blackboard framed diagrams.} arc $b$, and $c$ is the arc on 
the left-hand side of $b$.  The axioms are then derived from the Reidemeister 
moves, resulting in an algebraic structure which is an invariant of knot type. 
It is well-known that the knot quandle is a complete invariant of classical 
knot type up to reflection, though in general quandle-equivalent knots need 
not be ambient isotopic.

One way of strengthening the resulting invariant is to repeat the procedure 
used to derive the quandle definition from knot diagrams with \textit{semiarcs}
in place of arcs -- instead of dividing our oriented knot diagram into arcs by 
breaking the diagram at the undercrossing points, we now also break it at 
overcrossing 
points. These semiarcs are oriented edges in the underlying 4-valent graph of 
the knot diagram. If we think of the two inbound semiarcs as the inputs to 
binary operations, then at each crossing we have two output 
semiarcs; since there are two types of crossings, this yields four binary 
operations. Comparing the semiarc labels around the edge of the circle before
and after a minimal set of oriented Reidemeister moves gives us a set of 
axioms. The resulting algebraic structure is called a \textit{biquandle}; see 
\cite{CES}, \cite{FRC}, \cite{FJK}, and \cite{KR}, for example.

In this paper, we study finite biquandles using a matrix notation similar to 
the quandle matrix notation in \cite{HN}; these finite biquandles define 
invariants of virtual and classical knots via homomorphism counting much like 
finite groups and finite quandles. The paper is organized as follows: in 
section \ref{sec:biq}, we give the biquandle definition and prove some results 
about finite biquandles. In section \ref{sec:inv} we define the biquandle 
counting invariant and our biquandle matrix notation. In section 
\ref{sec:invex} we exhibit finite biquandles whose counting invariants 
distinguish the virtual trefoil and various Kishino knots from the unknot, 
and a biquandle whose counting invariant 
distinguishes a virtual knot from its reflection and from its reverse 
(sometimes called the \textit{knot inverse}). In section \ref{last} we 
describe our algorithms
for finding finite biquandles and computing the counting invariant, and we give
a classification of finite biquandles with up to 4 elements. A file containing 
the Maple code used to obtain these results can be obtained at the first 
author's website at \texttt{www.esotericka.org/quandles}. We wish to thank 
the referee for helpful suggestions and comments.

\section{\large\textbf{Biquandles}} \label{sec:biq}

We begin by recalling a definition from \cite{KR}.

\begin{definition} \textup{
A \textit{biquandle} is a set $B$ with four binary operations
$B\times B\to B$ denoted by 
\[(a,b) \mapsto a^b, \ a^{\bar b}, \ a_b,\quad \mathrm{and} \quad a_{\bar b}\] 
respectively, satisfying the following 20 axioms:}

\begin{list}{}{}
\item[\textup{(1.)}]{\textup{For every pair of elements $a,b\in B$, we have}
\[
\mathrm{(i)} \ a=a^{b{\bar{b_a}}}, \quad       
\mathrm{(ii)} \ b=b_{a{\bar{a^b}}}, \quad  
\mathrm{(iii)} \ a=a^{\bar{b}b_{\bar a}}, \quad
\mathrm{and} \quad
\mathrm{(iv)} \ b=b_{\bar{a}a^{\bar b}}.
\]}
\item[\textup{(2.)}]{\textup{ Given elements $a,b\in B$, there are elements $x,y\in B$ such that}
\[ 
\mathrm{(i)} \ x=a^{b_{\bar x}}, \quad
\mathrm{ (ii) } \ a=x^{\bar b},\quad
\mathrm{ (iii)} \   b=b_{{\bar x}a}, \quad
\mathrm{ (iv)} \  y=a^{\bar {b_y}},\quad
\mathrm{ (v) } \  a=y^b,   \quad \mathrm{and} \quad
\mathrm{ (vi)} \  b=b_{y{\bar a}}.
\]}
\item[\textup{(3.)}]{\textup{ For every triple $a,b,c \in B$ we have:}
\[
\mathrm{ (i)} \ a^{bc}=a^{c_bb^c}, \quad
\mathrm{(ii)} \ c_{ba} =c_{a^bb_a}, \quad 
\mathrm{(iii)} \ (b_a)^{c_{a^b}}=(b^c)_{a^{c_b}}, \]
\[
\mathrm{(iv)} \ a^{{\bar b}{\bar c}}=a^{{\bar c}_{\bar b}{\bar b}^{\bar c}}, 
\quad
\mathrm{(v)} \ c_{{\bar b}{\bar a}} =c_{{\bar a}^{\bar b}{\bar b}_{\bar a}},
\quad \mathrm{and}  \quad
\mathrm{(vi)} \ (b_{\bar a})^{{\bar c}_{{\bar a}^{\bar b}}} =(b^{\bar c})_{{\bar a}^{{\bar c}_{\bar b}}}.
\]}
\item[\textup{(4.)}]{\textup{Given an element $a\in B,$ there are elements 
$x,y\in B$ such that}
\[
\mathrm{ (i)} \   x=a_x, \quad 
\mathrm{ (ii)} \  a=x^a, \quad 
\mathrm{(iii)} \  y=a^{\bar y}, \quad \mathrm{and}  \quad 
\mathrm{(iv) } \  a=y_{\bar a}.
\]
}
\end{list}\textup{
The operations with bars are called \textit{left operations} and the operations
without bars are \textit{right operations.} The operations denoted by 
subscripts are \textit{lower operations} while those denoted by superscripts
are \textit{upper operations.}}
\end{definition}

This definition is obtained by dividing an oriented knot or link diagram at 
every crossing point to obtain a collection of \textit{semiarcs}. These 
semiarcs are the edges of the knot diagram considered as a 4-valent graph 
enhanced with crossing information. The two inbound semiarcs then operate on 
each other to yield the two outbound semiarcs, with different operations at 
positive and negative crossings as depicted in figure 1. The axioms are then 
transcriptions of a minimal set of oriented Reidemeister moves which are 
sufficient to generate any other oriented Reidemeister move. Axioms (1) and 
(2) come from the direct and reverse type II moves respectively, axiom (3) 
comes from the two type III moves with all positive and all negative crossings,
and axiom (4) comes from the type I move. 

\begin{figure}
\[\includegraphics{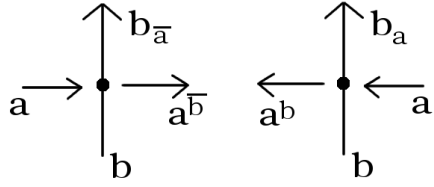}\]
\caption{Biquandle operations at crossings}
\end{figure}

\begin{remark}\textup{ \label{switch}
There is an alternate notation (see \cite{CES} and \cite{FJK}) which is both 
more general and simpler in some ways than the biquandle notation we use in 
this paper. Specifically, a \textit{switch} on $X$ is a map 
$S:X\times X\to X\times X$ of the form $S(a,b)=(b_a,a^b)$ which satisfies
the \textit{Yang-Baxter equation}
\[(S\times \mathrm{Id})(\mathrm{Id}\times S)(S\times \mathrm{Id}) =
(\mathrm{Id}\times S)(S\times \mathrm{Id})(\mathrm{Id}\times S).
\] The direct Reidemeister II move then requires that $S$ be invertible, 
with the inverse given by $S^{-1}(a,b)=(b^{\bar a}, a_{\bar b})$.
A biquandle is then a switch which 
satisfies the extra conditions coming from the type I and reverse type II 
moves. For the purpose of finding finite biquandles by filling in their
matrices, as we will describe later, the notation in definition 1 is perhaps 
more useful.
}\end{remark}

\begin{proposition}
Let $B$ be a biquandle and define the \textit{obverse} of $B$, $Obv(B)$, to be 
the algebra obtained from $B$ by interchanging the right operations with the
left ones, and let the \textit{flip} of $B$, $Flip(B)$, be the algebra
obtained by interchanging the upper operations with the lower operations. Then
$\mathrm{Obv}(B)$, $\mathrm{Flip}(B)$ and $\mathrm{Obv}(\mathrm{Flip}(B))
=\mathrm{Flip}(\mathrm{Obv}(B))$ are also biquandles.
\end{proposition}

\begin{proof}
Here we only need to observe that the set of axioms is symmetrical under the
operation of interchanging the right and left operations and upper and lower 
operations. Inspection shows that this hold for axioms (1), (3) and (4).
To see that axiom (2) also satisfies this symmetry,
consider the following equivalent reformulation of axiom (2):
\begin{list}{}{}
\item[\textup{(2$'$.)}]{\textup{ Given elements $a,b\in B$, there are 
elements $x,y,c,d\in B$ such that}
\[ 
\mathrm{(i)} \ c=b_{\bar x}, \quad
\mathrm{(ii)} \ x=a^c, \quad
\mathrm{ (iii) } \ a=x^{\bar b},\quad
\mathrm{ (iv)} \   b=c_a, \quad \mathrm{and}\]\[
\mathrm{(v)} \  d=b_y, \quad
\mathrm{ (vi)} \  y=a^{\bar d},\quad
\mathrm{ (vi) } \  a=y^b,   \quad 
\mathrm{ (vii)} \  b=d_{\bar a}.
\]}
\end{list}
\end{proof}

We will also need the following definition (see \cite{J} or \cite{FR}):

\begin{definition}\textup{
A \textit{quandle} is a set $Q$ with a binary operation $(a,b)\mapsto a^b$
such that
\newcounter{qax}
\begin{list}{(\roman{qax})}{\usecounter{qax}}
\item{For all $a\in Q$, $a^a=a$,}
\item{For all $a,b \in Q$, there is a unique $c\in Q$ such that $a=c^b$, and}
\item{For all $a,b,c\in Q$, $a^{bc}=a^{cb^c}$.}
\end{list}
The uniqueness in axiom (ii) implies that the map $f_b:Q\to Q$ given by
$f_b(a)=a^b$ is a bijection; the inverse map is denoted 
$f^{-1}_b(a)=a^{\bar b}$, and it is an easy exercise to show that $Q$ forms 
a quandle under the operation $(a,b)\to a^{\bar b}$, called the \textit{dual}
of $Q$. Axiom (ii) may then be reformulated as
\begin{list}{}{}
\item[(ii$'$)]{For every $a,b\in Q$ we have $a^{b\bar{b}}=a^{\bar{b}b}=a.$}
\end{list}
}\end{definition}

As expected, we have the following:

\begin{definition}\textup{
A map $\phi:B_1\to B_2$ is a \textit{biquandle homomorphism} if it preserves 
all four biquandle operations, that is 
\[ \phi(a^b)=\phi(a)^{\phi(b)}, \phi(a^{\bar b})=\phi(a)^{\overline{\phi(b)}}, 
 \phi(a_b)=\phi(a)_{\phi(b)}, \quad \mathrm{and} \quad
\phi(a_{\bar b})=\phi(a)_{\overline{\phi(b)}}.\]
A bijective biquandle homomorphism is a \textit{biquandle isomorphism.}}
\end{definition}

The knot biquandle of the obverse of a knot is the obverse of the biquandle of 
the original knot. The knot biquandle of the flip of a knot, that is, of the 
knot viewed from the other side of its supporting surface, is the flip of the 
biquandle of the original knot. If a knot is classical, then flipping is an 
ambient isotopy and the resulting biquandle must be isomorphic to its flip. 
If a virtual knot is non-classical, however, the flip of its biquandle may not 
be isomorphic to the original knot biquandle.

\begin{definition}\textup{
A biquandle $B$ is \textit{self-obverse} if $B$ is isomorphic to 
$\mathrm{Obv}(B)$. $B$ is \textit{self-Flip} if $B$ is isomorphic to 
$\mathrm{Flip}(B)$.}
\end{definition}

\begin{proposition}
Let $B$ be a biquandle. Then the automorphism groups $\mathrm{Aut}(B)$,
$\mathrm{Aut}(\mathrm{Obv}(B))$, $\mathrm{Aut}(\mathrm{Flip}(B))$ and
$\mathrm{Aut}(\mathrm{Flip}(\mathrm{Obv}(B)))=
\mathrm{Aut}(\mathrm{Obv}(\mathrm{Flip}(B)))$ are all isomorphic.
\end{proposition}

\begin{proof}
Let $\phi:B\to B$ be an automorphism. Then for all $a,b\in B$, we have
$\phi(a^b)=\phi(a)^{\phi(b)},$ $\phi(a^{\bar b})=\phi(a)^{\bar{\phi(b)}},$ 
$\phi(a_b)=\phi(a)_{\phi(b)}$, and $\phi(a_{\bar b})=\phi(a)_{\bar{\phi(b)}}$.
Then when viewed as a map from $\mathrm{Obv}(B)$ to itself, switching
the right and left operations both before and after applying the map $\phi$
gives us back the same set of equalities, and 
$\phi\in \mathrm{Aut}(\mathrm{Obv}(B))$. Similarly for the other cases.
\end{proof}

To define biquandle presentations in terms of generators and relations, for a 
finite set $X$ we start by defining the set of \textit{biquandle words in $X$},
$BW(X)$, to be the set which includes all elements of $X$ together with 
all finite strings of the forms
\[ a^b, \ a^{\bar b}, \ a_b, \ a_{\bar b}, \ X(a,b), \ Y(a,b), \ Fx(a) \ 
\mathrm{and} \ Fy(a)\]
where $a,b \in BW(X)$.
Then the \textit{free biquandle on $X$} is the set of equivalence classes in
$BW(X)$ under the equivalence relation generated by the relations required by 
the biquandle axioms above, e.g. $a^{bc}\sim a^{c_bb^c}$, 
$X(a,b)\sim a^{b_{\bar{X(a,b)}}} $, $Fx(a)\sim a_{Fx(a)}$, etc. 
A \textit{finitely 
presented biquandle} is then the quotient of the free biquandle on a finite 
set $X$ by the equivalence relation generated by a finite list of equivalences 
of biquandle words, which we may call \textit{explicit relations} to 
distinguish them from the \text{implicit relations} required by the biquandle 
axioms. See \cite{KM} for more on the universal algebra of biquandles.

For any virtual knot diagram $K$, we obtain a presentation of the \textit{knot
biquandle} by assigning a distinct generator to every semiarc and obtaining
a pair of relations at every crossing according to the diagrams in figure 1.
The knot biquandle is then quotient of the free biquandle on the set of 
semi-arcs in the virtual knot diagram by the set of relations at each 
crossing; it is an invariant of virtual isotopy by construction. 

\begin{definition}\textup{
A biquandle relation is \textit{short} if it is of the form $a=b^c$, 
$a=b^{\bar c}$, $a=b_c$ or $a=b_{\bar c}.$ In each of these, $a$ is the
\textit{input}, $b$ is the \textit{operator} and $c$ is the \textit{output}.
A biquandle presentation is
\textit{knotlike} if 
\newcounter{kl}
\begin{list}{(\roman{kl})}{\usecounter{kl}}
\item{Every explicit relation is short,}
\item{Every generator appears exactly once each as input, operator and output,}
\item{The explicit relations come in pairs of the form ($a^b=c, b_a=d$) or 
($a^{\bar b}=c, b_{\bar a}=d$). That is, if $a$ operates on input $b$ in one 
relation, then $b$ operates on input $a$ in another relation where one is an 
upper and the other a lower operation, and both are right- or both are left- 
operations.}
\end{list}
A biquandle is \textit{knotlike} if it has a knotlike presentation.
}\end{definition}

\begin{remark}\textup{
Condition (iii) says that the relations are \textit{switch relations}
$S(a,b)=(x,y)$. 
}\end{remark}

In particular, if we number the semiarcs in a knot biquandle sequentially 
following the orientation of the knot, then every relation is of the form
\[ i_y=i+1,\quad i_{\bar y}=i+1,  \quad i^y=i+1, \quad \mathrm{or} \quad 
i^{\bar y}=i+1   \]
and we can unambiguously specify a knot biquandle as a vector whose $i$th
entry is the operation in the explicit relation with input $i$. This vector 
notation can easily be translated into a Gauss code (\cite{K}) or an SOKQ 
presentation (\cite{N2}).

Every virtual knot or link biquandle is knotlike. Unlike the knot quandle,
given a knotlike biquandle presentation we can reconstruct the virtual knot or 
link it comes from up to strictly virtual moves. This is consistent with the 
following conjecture (see \cite{FJK}):

\begin{conjecture}
The knot biquandle is complete invariant of virtual link type. That is, if 
two knotlike biquandle presentations are related by Tietze moves, then the 
resulting virtual link diagrams are related by virtual isotopy moves.
\end{conjecture}

As with other algebraic invariants of knots and links, direct comparison of 
isomorphism types of biquandles given by presentations is generally difficult. 
Thus, we seek invariants of biquandle isomorphism type, which are naturally
also knot invariants. One such invariant is the counting invariant associated 
to a finite biquandle, described in the following section.

\section{\large \textbf{Finite biquandles and the counting invariant}}
\label{sec:inv}

Let $T$ be a finite biquandle and $K$ a finitely presented biquandle, 
e.g. a knot biquandle. As with groups and quandles, we have the following
theorem (see also \cite{FJK}):

\begin{theorem}
The cardinality $|\mathrm{Hom}(K,T)|$ of the set of homomorphisms from $K$ to 
$T$ is an invariant of biquandle isomorphism type and hence an invariant of 
virtual isotopy.
\end{theorem} 

\begin{proof}
If $\phi:K\to K'$ is an isomorphism then for every $f\in \mathrm{Hom}(K,T)$ 
we have $f\phi^{-1}\in \mathrm{Hom}(K',T)$ and for every 
$g\in \mathrm{Hom}(K',T)$ we have $g\phi \in \mathrm{Hom}(K,T)$. 
Hence we have both $|\mathrm{Hom}(K,T)|\le |\mathrm{Hom}(K',T)|$ and
$|\mathrm{Hom}(K',T)|\le |\mathrm{Hom}(K,T)|$, so 
$|\mathrm{Hom}(K,T)|= |\mathrm{Hom}(K',T)|$.
\end{proof}

Thus, for every finite biquandle $T$, the number of homomorphisms into $T$ is
an invariant of virtual (and hence classical) knot type. In order to 
find and evaluate such invariants, then, we need a way of representing finite 
biquandles. In \cite{HN}, finite quandles $Q=\{x_1,\dots, x_n\}$ of order $n$ 
are represented as $n\times n$ matrices with $M_{ij}=k$ where $x_k=x_i^{x_j}$.
To represent finite biquandles as matrices, we use the following block matrix
notation.

\begin{definition}\textup{
Let $B=\{x_1,\dots, x_n\}$ be a finite biquandle. Then the \textit{matrix of B}
is the block matrix}
\[
M_B=\left[\begin{array}{l|l} M^1 & M^2 \\ \hline M^3 & M^4 
\end{array}
\right]
\qquad \mathrm{where} \qquad 
M^l_{ij} =k_l,\quad 
x_{k_l}=\left\{\begin{array}{ll}
x_i^{\bar{x_j}}, & l=1, \\
x_i^{x_j}, & l=2, \\
(x_i)_{\bar{x_j}}, & l=3, \\
(x_i)_{x_j}, & l=4.
\end{array}\right.\]
\textup{
The matrices $M^i$ will be called the \textit{block submatrices} or just the
\textit{submatrices} of $B$.}
\end{definition}

\begin{remark}\textup{
In light of remark \ref{switch}, we can without loss of information 
drop the two left matrices to obtain a \textit{switch matrix}
\[
M_S=\left[ M^2 | M^4 \right].
\]
In particular, from a switch matrix we can recover the biquandle matrix.
}\end{remark}

\begin{example} \textup{
For any positive integer $n$ let $T_n$ be the matrix such that the $i$th row
of $T_n$ has every entry equal to $i$. Then the matrix}
\[BT_n =\left[\begin{array}{c|c} T_n & T_n \\ \hline T_n & T_n 
\end{array}\right]\]
\textup{is a finite biquandle, called the \textit{trivial biquandle} of 
order $n$. For instance, }
\[
BT_2 = \left[\begin{array}{cc|cc} 1 & 1 & 1 & 1 \\ 2 & 2 & 2 & 2 \\ \hline 
1 & 1 & 1 & 1 \\ 2 & 2 & 2 & 2 
\end{array}\right],\quad
BT_3 = \left[\begin{array}{ccc|ccc} 
1 & 1 & 1 & 1 & 1 & 1 \\ 2 & 2 & 2 & 2 & 2 & 2 \\ 3 & 3 & 3 & 3 & 3 & 3 
\\ \hline 
1 & 1 & 1 & 1 & 1 & 1 \\ 2 & 2 & 2 & 2 & 2 & 2 \\ 3 & 3 & 3 & 3 & 3 & 3 
\end{array}\right], \quad \mathrm{etc.}
\]
\end{example}

Even though the biquandle axioms, unlike the quandle axioms, do not require
uniqueness of right inverses under the four biquandle operations, we have

\begin{proposition} \label{fb1}
If $B$ is a finite biquandle, then the right inverses under each of the four
actions are unique. 
\end{proposition}

\begin{proof}
Axiom (2) part (v) says that for every $a,b\in B$ there is a $c$ such that 
$a=c^b$.
That is, for a given element $b\in B$, every element of $B$ must appear at 
least once in the column of $b$ in the $M^2$ block of the matrix of $B$. Since 
there are only $n$ available positions and $n$ elements which must appear,
every element appears exactly once. Similarly for $M^1$.

To see uniqueness for $M^4$, note that axiom 2 part (iii) is equivalent to
``for all $a,b\in B$ there exist $c,x\in B$ such that $b=c_a$ and 
$c=b_{\bar x}$.'' In particular, every element $b$ must appear in some
row $c$ in column $a$, and as before we have uniqueness. Similarly for $M^3$.
\end{proof}

Thus, we have one easily checked condition for a block matrix to be the matrix
of a finite biquandle, namely that each block must have columns which are
permutations. The condition arising from the type I move is also easy to
check visually:

\begin{proposition} \label{fb2}
If $M$ is one of the four submatrices of a finite biquandle, then every row 
of $M$ has exactly one entry equal to its column number in $M$. 
\end{proposition}

\begin{proof}
Axiom (4) implies that every row contains at least one entry which equals its
column number. Now suppose row $i$ has $M_{ij}=j$ and $M_{ik}=k$. Then since
every column is a permutation, no entry appears more than once in a column, so
no other row has a $j$ in column $j$ or a $k$ in column $k$. Then of the 
remaining $n-1$ rows, at most $n-2$ have an entry equal to its column number,
and $M$ is not a biquandle submatrix.
\end{proof}

\begin{proposition}
Let $B$ be a finite biquandle. If the lower submatrices of $B$ are both trivial
matrices $T_n$, i.e. if $x_y=x$ and $x_{\bar y}=x$ for all $x,y \in B$, then
the upper submatrices of $B$ form a quandle matrix and its dual.
\end{proposition}

\begin{proof}
If the lower operations are trivial, then the biquandle axioms reduce to
\begin{list}{}{}
\item[1.]{For every $a,b\in B$ we have $a=a^{b\bar{b}}$ and $a=a^{\bar{b} b}.$}
\item[2.]{For all $a,b\in B$ there is an $x$ such that $a=x^{\bar b}$ and
 a $y$ such that $a=y^b$.}
\item[3.]{For every triple $a,b,c\in B$ we have $a^{bc}=a^{cb^c}$ and 
$a^{\bar{b} \bar{c}} =a^{\bar{c} \bar{b}^{\bar{c}}}$.}
\item[4.]{For every $a$ in $B$, we have $a=a^a$ and $a=a^{\bar a}.$}
\end{list}
But these are just the quandle axioms; (1) and (2) together give us quandle 
axiom (ii), while (3) is quandle axiom (iii) for a quandle and its dual, and
(4) is quandle axiom (i) for a quandle and its dual.
\end{proof}

\begin{corollary}
Let $Q$ be a finite quandle. Then there is a corresponding finite biquandle
with matrix
\[B=\left[\begin{array}{c|c} M_{\bar Q} & M_Q \\ \hline T_n  & T_n
\end{array}\right],\]
where $M_Q$ is the matrix of $Q$, and $M_{\bar Q}$ is the matrix of the dual 
of $Q$, and $T_n$ is the matrix of the trivial quandle of order $n$. Such a
biquandle is self-obverse iff $Q$ is self-dual, and self-flip iff $Q$ is 
trivial.
\end{corollary}

A biquandle with both lower (or both upper) operations trivial is really a 
quandle in disguise, which we call a \textit{quandle biquandle} or a 
\textit{qbiq}. Thus, to find new knot invariants using finite biquandles, we 
seek biquandles in which the lower (or upper) operations are non-trivial; we 
call such biquandles \textit{non-quandle biquandles} or \textit{non qbiqs}.

In the next section we give a few examples to illustrate the usefulness of
biquandle counting invariants.

\section{\large\textbf{Counting invariant examples}}\label{sec:invex}

To compute the counting invariant we first obtain a presentation of the
knot biquandle for the virtual knot or link in question. We represent such
a presentation with a \textit{biquandle presentation matrix} in which each 
short-form relation fills in an entry in one of the biquandle submatrices, and
the remaining entries are zero. Thus, an $n$-crossing knot or link has a 
$4n\times 4n$ biquandle presentation matrix with exactly $2n$ non-zero entries.

\begin{example}\label{un}
\textup{The biquandle of the unknot is the free biquandle on one generator. 
The zero-crossing diagram of the unknot determines a biquandle 
presentation with a single generator and an empty list of explicit relations, 
which does not yield a useful presentation matrix. The below one-crossing 
unknot diagram has biquandle presentation $\langle 1,2 \ | \ 1^{\bar 2}=1, 
2_{\bar 1}=1 \rangle$, presentation vector $[{}^{\bar 2},{}_{\bar 1}]$ and 
presentation matrix}
\[
\raisebox{-0.5in}{\includegraphics{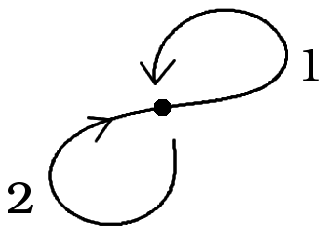}} \quad \quad
U=\left[\begin{array}{cc|cc}
0 & 2 & 0 & 0 \\
0 & 0 & 0 & 0 \\ \hline
0 & 0 & 0 & 0 \\
1 & 0 & 0 & 0 
\end{array}
\right].
\]
\end{example}

\begin{example}\label{vt} 
\textup{The virtual trefoil from \cite{K} has biquandle presentation
$\langle 1,2,3,4 \ | \ 1_{\bar 3}=2, 2_{\bar 4}=3, 3^{\bar 1}=4, 
4^{\bar 2}=1\rangle$, which has presentation vector $[{}_{\bar 3},\ 
{}_{\bar 4},\ {}^{\bar 1},\ {}^{\bar 2}]$ and presentation matrix}
\[
\raisebox{-0.5in}{\includegraphics{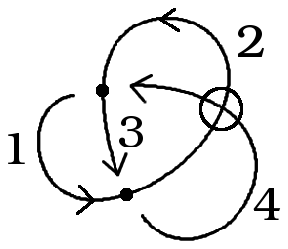}} \quad \quad
VT=\left[\begin{array}{cccc|cccc}
0 & 0 & 0 & 0 & 0 & 0 & 0 & 0 \\
0 & 0 & 0 & 0 & 0 & 0 & 0 & 0 \\
4 & 0 & 0 & 0 & 0 & 0 & 0 & 0 \\
0 & 1 & 0 & 0 & 0 & 0 & 0 & 0 \\ \hline
0 & 0 & 2 & 0 & 0 & 0 & 0 & 0 \\
0 & 0 & 0 & 3 & 0 & 0 & 0 & 0 \\
0 & 0 & 0 & 0 & 0 & 0 & 0 & 0 \\
0 & 0 & 0 & 0 & 0 & 0 & 0 & 0 
\end{array}
\right].
\]
\end{example}

The counting invariant for a non-quandle biquandle can give us information not 
present in the knot quandle. For instance, the virtual trefoil in example 
\ref{vt} has trivial knot quandle (though it has nontrivial Jones polynomial); 
however, the biquandle counting invariant with target biquandle
\[
T=\left[
\begin{array}{ccc|ccc}
2 & 1 & 3 & 3 & 2 & 1 \\
1 & 3 & 2 & 2 & 1 & 3 \\
3 & 2 & 1 & 1 & 3 & 2 \\
\hline
2 & 2 & 2 & 3 & 3 & 3 \\
3 & 3 & 3 & 1 & 1 & 1 \\
1 & 1 & 1 & 2 & 2 & 2 \\ 
\end{array}
\right]
\]
distinguishes the virtual trefoil from the unknot, as we have
$|\mathrm{Hom}(U,T)|=3\ne 0=|\mathrm{Hom}(VT,T)|$ where $U$ and $VT$ are
the biquandles of the unknot and the virtual trefoil as described above.

\begin{example}
\textup{The pictured Kishino knot has trivial (upper) quandle (\cite{KS}). 
Its biquandle
has presentation \[\langle 1,2,3,4,5,6,7,8 \ | \ 1_7=2, 2^{\bar 8}=3, 3^5=4,
4_{\bar 6}=5, 5_3=6, 6^{\bar 4}=7, 7^1=8, 8_{\bar 2}=1 \rangle\]
with presentation vector $[{}_7, {}^{\bar 8}, {}^5,
{}_{\bar 6}, {}_3, {}^{\bar 4}, {}^1, {}_{\bar 2}]$ and presentation matrix}
\[
\raisebox{-0.5in}{\includegraphics{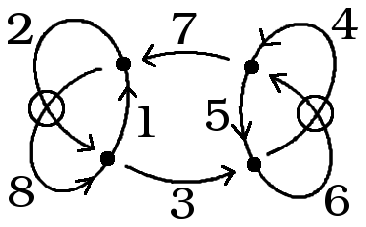}} \quad \quad
K=\left[\begin{array}{cccccccc|cccccccc}
0 & 0 & 0 & 0 & 0 & 0 & 0 & 0 &   0 & 0 & 0 & 0 & 0 & 0 & 0 & 0  \\
0 & 0 & 0 & 0 & 0 & 0 & 0 & 3 &   0 & 0 & 0 & 0 & 0 & 0 & 0 & 0  \\
0 & 0 & 0 & 0 & 0 & 0 & 0 & 0 &   0 & 0 & 0 & 0 & 4 & 0 & 0 & 0  \\
0 & 0 & 0 & 0 & 0 & 0 & 0 & 0 &   0 & 0 & 0 & 0 & 0 & 0 & 0 & 0  \\
0 & 0 & 0 & 0 & 0 & 0 & 0 & 0 &   0 & 0 & 0 & 0 & 0 & 0 & 0 & 0  \\
0 & 0 & 0 & 7 & 0 & 0 & 0 & 0 &   0 & 0 & 0 & 0 & 0 & 0 & 0 & 0  \\
0 & 0 & 0 & 0 & 0 & 0 & 0 & 0 &   8 & 0 & 0 & 0 & 0 & 0 & 0 & 0  \\
0 & 0 & 0 & 0 & 0 & 0 & 0 & 0 &   0 & 0 & 0 & 0 & 0 & 0 & 0 & 0  \\ \hline
0 & 0 & 0 & 0 & 0 & 0 & 0 & 0 &   0 & 0 & 0 & 0 & 0 & 0 & 2 & 0  \\
0 & 0 & 0 & 0 & 0 & 0 & 0 & 0 &   0 & 0 & 0 & 0 & 0 & 0 & 0 & 0  \\
0 & 0 & 0 & 0 & 0 & 0 & 0 & 0 &   0 & 0 & 0 & 0 & 0 & 0 & 0 & 0  \\
0 & 0 & 0 & 0 & 0 & 5 & 0 & 0 &   0 & 0 & 0 & 0 & 0 & 0 & 0 & 0  \\
0 & 0 & 0 & 0 & 0 & 0 & 0 & 0 &   0 & 0 & 6 & 0 & 0 & 0 & 0 & 0  \\
0 & 0 & 0 & 0 & 0 & 0 & 0 & 0 &   0 & 0 & 0 & 0 & 0 & 0 & 0 & 0  \\
0 & 0 & 0 & 0 & 0 & 0 & 0 & 0 &   0 & 0 & 0 & 0 & 0 & 0 & 0 & 0  \\
0 & 1 & 0 & 0 & 0 & 0 & 0 & 0 &   0 & 0 & 0 & 0 & 0 & 0 & 0 & 0 
\end{array}
\right].\] 
\textup{This virtual knot is distinguished from the unknot 
by the biquandle counting invariant with target biquandle
\[
T_2=\left[
\begin{array}{cccc|cccc}
3 & 1 & 2 & 4 & 4 & 1 & 3 & 2 \\
2 & 4 & 3 & 1 & 2 & 3 & 1 & 4 \\
1 & 3 & 4 & 2 & 3 & 2 & 4 & 1 \\
4 & 2 & 1 & 3 & 1 & 4 & 2 & 3 \\ 
\hline
4 & 1 & 3 & 2 & 3 & 1 & 2 & 4 \\
2 & 3 & 1 & 4 & 2 & 4 & 3 & 1 \\
3 & 2 & 4 & 1 & 1 & 3 & 4 & 2 \\
1 & 4 & 2 & 3 & 4 & 2 & 1 & 3  
\end{array}
\right],\]
as $|\mathrm{Hom}(K,T_2)|=16 \ne |\mathrm{Hom}(U,T_2)|=4$.
}
\end{example}

\begin{example}
\textup{The Kishino knots below both have a biquandle counting invariant 
value of 16 with the finite biquandle $T_4$ below, 
while $|\mathrm{Hom}(U,T_4)|=4$: }

\begin{center}
\parbox{2in}{\includegraphics{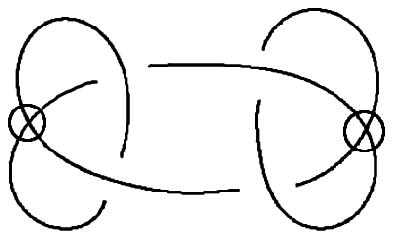}  

\includegraphics{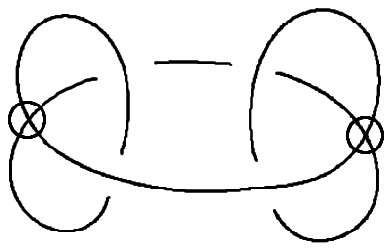}}
\parbox{3in}{\quad \quad $
T_4=\left[
\begin{array}{cccc|cccc}
1 & 4 & 2 & 3 & 1 & 3 & 4 & 2 \\
2 & 3 & 1 & 4 & 3 & 1 & 2 & 4 \\
4 & 1 & 3 & 2 & 2 & 4 & 3 & 1 \\
3 & 2 & 4 & 1 & 4 & 2 & 1 & 3 \\ 
\hline
1 & 3 & 4 & 2 & 1 & 4 & 2 & 3 \\
3 & 1 & 2 & 4 & 2 & 3 & 1 & 4 \\
2 & 4 & 3 & 1 & 4 & 1 & 3 & 2 \\
4 & 2 & 1 & 3 & 3 & 2 & 4 & 1  
\end{array}
\right].$}
\end{center}

\end{example}

\begin{example}\textup{This Kishino knot has previously been shown 
to be non-trivial using quaternionic biquandles \cite{BF}; however, our
method detects its non-triviality with the four-element biquandle
$T_4$ from the previous example, again with $|\mathrm{Hom}(K,T_4)|=16$. 
For completeness, we provide a list of the homomorphisms from the knot 
biquandle to $T_4$ as computed by our program. We would like to thank the 
referee for suggesting this example.}

\medskip

\textup{
\raisebox{-0.5in}{\includegraphics{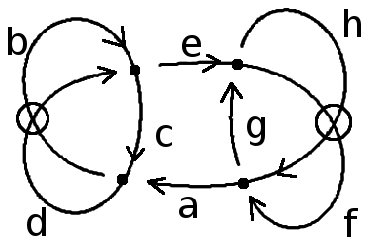}}
\quad
\begin{tabular}{|cccccccc|} \hline
a & b & c & d & e & f & g & h\\ \hline
1 & 1 & 1 & 1 & 1 & 1 & 1 & 1 \\
1 & 1 & 1 & 1 & 1 & 2 & 3 & 2 \\
1 & 2 & 3 & 2 & 1 & 1 & 1 & 1 \\
1 & 2 & 3 & 2 & 1 & 2 & 3 & 2 \\
2 & 1 & 3 & 4 & 4 & 3 & 2 & 1 \\
2 & 1 & 3 & 4 & 4 & 4 & 3 & 1 \\
2 & 2 & 1 & 3 & 4 & 3 & 1 & 2 \\
2 & 2 & 1 & 3 & 4 & 4 & 3 & 1 \\ \hline
\end{tabular}
\quad 
\begin{tabular}{|cccccccc|} \hline
a & b & c & d & e & f & g & h\\ \hline
3 & 3 & 3 & 3 & 3 & 3 & 3 & 3 \\
3 & 3 & 3 & 3 & 3 & 4 & 1 & 4 \\
3 & 4 & 1 & 4 & 3 & 3 & 3 & 3 \\
3 & 4 & 1 & 4 & 3 & 4 & 1 & 4 \\
4 & 3 & 1 & 2 & 2 & 1 & 3 & 4 \\
4 & 3 & 1 & 2 & 2 & 2 & 1 & 3 \\
4 & 4 & 3 & 1 & 2 & 1 & 3 & 4 \\
4 & 4 & 3 & 1 & 2 & 2 & 1 & 3 \\ \hline
\end{tabular}
}
\end{example}

\bigskip

Note that, unlike quandles and groups, for a given pair of biquandles there 
may or may not be a homomorphism between the pair; thus a biquandle counting 
invariant can take any non-negative integer value, including zero. Indeed,
we have

\begin{proposition}
Let $T$ be a finite biquandle and let $N(T)$ be the number of elements
of $T$ which are idempotent in all four operations. Then for any finitely
presented biquandle $B$, the counting invariant with target $T$ satisfies
\[|\mathrm{Hom}(B,T)|\ge N(T).\]
\end{proposition}

\begin{proof}
If $a\in B$ is idempotent in all four operations, that is, if
\[ a^{\bar a}=a, a^a=a, a_{\bar a}=a \ \mathrm{and} \ a_a=a,\]
then for any biquandle $B$ with $n$ generators, the constant map sending all
generators to $a$ satisfies any possible list of relations and hence is always
a biquandle homomorphism.
\end{proof}

Attempting to compute a complete list of finite biquandles even for small
order quickly becomes computationally very resource-hungry. To find suitable 
finite target biquandles for use in counting invariants, then, it can be 
helpful to look for \textit{biquandle completions}.

\begin{definition}\textup{
Let $M_B$ be a biquandle presentation matrix. A \textit{biquandle completion} 
of $M_B$ is a finite biquandle whose matrix is obtained by filling in the 
zeroes of $M_B$.}
\end{definition}

Since filling in a zero amounts to adding a new explicit short relation,
a completion is actually a finite quotient of the original biquandle. In 
particular, a completion of the presentation matrix of a knotlike biquandle 
$KB$ is a finite quotient biquandle $T$ onto which there is a surjective 
homomorphism. More generally, for any labelling of the semi-arcs in a
knot diagram with numbers $\{1,2,\dots, n\}$ such that no two crossings 
have the same input labels, we can form the resulting biquandle presentation 
matrix and look for completions. 

\begin{definition}\textup{
Let $B$ be a biquandle. If there is an element $a\in B$ such that every 
element of $B$ is equivalent to a word starting with $a$ (that is, a word with 
leftmost generator $a$), we say $B$ is \textit{connected}.
}\end{definition}

As with quandles, knot biquandles are connected, though link biquandles
in general are not. Hence, for finding new knot invariants, we are primarily 
concerned with finding finite connected non-quandle biquandles. Of particular 
interest are finite non-quandle biquandles which are not self-flip and not
self-obverse. As with quandles and groups, two isomorphic biquandles define 
the the same knot invariant; however, the flip or the obverse of a non 
self-flip or non self-obverse biquandle generally defines a different 
invariant from that defined by the original.

Our results in section \ref{last} indicate that all biquandles with order
$\le 4$ are self-obverse. Non self-obverse biquandles can be used to 
distinguish some virtual knots from their obverses -- if 
$|\mathrm{Hom}(K,T)| \ne |\mathrm{Hom}(Obv(K),T)|$ or equivalently if
$|\mathrm{Hom}(K,T)| \ne |\mathrm{Hom}(K,Obv(T))|$, then $K$ is not
isotopic to its obverse.

\begin{example} \textup{
The biquandle $T_5$ below is not self-obverse. Its counting invariant
distinguishes the virtual knot below from its obverse, with
$|\mathrm{Hom}(K_2,T_5)|=5 \ne 1 =|\mathrm{Hom}(Obv(K_2),T_5)|$. This 
virtual knot is similar to the virtual knots in figure 3 of \cite{JS}. 
Indeed, this same biquandle distinguishes $K_2$ from its \textit{inverse},
that is, the virtual knot obtained from $K_2$ by reversing its orientation.}
\[ K_2=
\raisebox{-0.5in}{\includegraphics{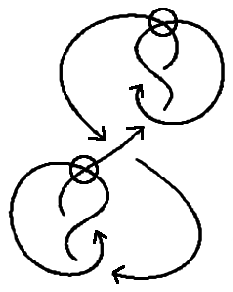}} \quad \quad
T_5=
\left[\begin{array}{ccccc|ccccc}
5 & 3 & 1 & 4 & 2 &  4 & 3 & 2 & 1 & 5 \\
4 & 2 & 5 & 3 & 1 &  3 & 2 & 1 & 5 & 4 \\
3 & 1 & 4 & 2 & 5 &  2 & 1 & 5 & 4 & 3 \\
2 & 5 & 3 & 1 & 4 &  1 & 5 & 4 & 3 & 2 \\
1 & 4 & 2 & 5 & 3 &  5 & 4 & 3 & 2 & 1 \\ \hline
5 & 5 & 5 & 5 & 5 &  4 & 4 & 4 & 4 & 4 \\
2 & 2 & 2 & 2 & 2 &  2 & 2 & 2 & 2 & 2 \\ 
4 & 4 & 4 & 4 & 4 &  5 & 5 & 5 & 5 & 5 \\
1 & 1 & 1 & 1 & 1 &  3 & 3 & 3 & 3 & 3 \\
3 & 3 & 3 & 3 & 3 &  1 & 1 & 1 & 1 & 1 \\
\end{array}\right].
\]
\end{example}

In the last section we describe an algorithm for finding finite biquandles,
including completions, and give a classification of all finite biquandles 
with up to 4 elements.

\section{\large\textbf{Computational Results}}\label{last}

A file called \texttt{biquandles-maple.txt} containing Maple code implementing 
the following algorithm is available for download at 
\texttt{www.esotericka.org/quandles}; we have chosen not to include it here
as it is rather lengthy. 

Our code represents a finite biquandle internally as a vector of four matrices 
rather than one large block matrix. Our first program, 
\texttt{biqtest}, checks a list of four $n\times n$ matrices for the biquandle
axioms, returning ``true'' if the list represents a finite biquandle and
``false'' if any axiom is not satisfied.

Much of the work in finding finite biquandles is done by \texttt{biqfill}, 
which takes as input a biquandle pattern consisting of four $n\times n$ 
matrices with 
entries in $\{0,1,\dots, n\}$. An entry of zero counts as a blank; our matrix 
notation lets the entries act both as biquandle elements and as submatrix 
row and column numbers. Overall, \texttt{biqfill} is a loop controlled by two 
variables, a ``changed'' counter and a ``contradiction'' counter. The program 
systematically checks a working biquandle pattern for each of the axioms. 
When the contradiction marker is set, the program skips the remaining checks, 
exits the loop and reports ``false''. Setting the ``changed'' counter to
``true'' tells the program to repeat the loop, propagating the changed value 
through the matrix. Biquandle axioms (1) and (3) each assert that various 
biquandle words corresponding to entries in the matrices should be equal; if 
enough entries are non-zero to determine the positions of the two entries in 
question, the entries are compared. If the two entries are equal, the 
algorithm moves on; if one is zero and the other non-zero, the zero entry is 
changed to match the non-zero entry and the ``changed'' counter is set
to true; if both entries are non-zero and not equal, then the ``contradiction''
counter is set to true. Axiom (2) makes assertions about the right inverses
of the various biquandle operations; \texttt{biqfill} tests for these, again
setting the contradiction counter if the pattern violates one of the axioms
and filling in zeroes when possible. For axiom (4), the contradiction counter 
is set to ``true'' if any row has no zeroes and no entry equal to its column 
number or if any row has more than one entry equal to its column number.
When a row has one entry equal to its column number, the corresponding entry 
in the upper or lower matrix is checked and filled in or the contradiction
counter is set as appropriate. 

The equations in axiom (3) come from the oriented versions of the third 
Reidemeister move with all three crossings of the same sign. To speed up the 
search algorithm, we include the equations arising from the other oriented 
Reidemeister III moves, which are consequences of axioms (1)-(3). Finally, 
\texttt{biqfill} checks each column, looking for any which have only
one zero; such zeroes are filled in or contradictions detected with the
program \texttt{avail}, which finds the smallest available entry for the
specified positions which does not contradict propositions \ref{fb1} and
\ref{fb2} or returns ``false'' if no such entry exists. If the 
``contradiction'' counter evaluates to ``true'' at any point, the loop is 
exited and \texttt{biqfill} returns ``false.'' At the end of the 
loop, if the ``changed'' counter evaluates to ``true'' then the loop is 
repeated; if not, the program reports the new pattern and exits. 

The program \texttt{biqlist} takes a starting biquandle pattern and starts a 
working list. The program removes the first biquandle pattern from the
working list, then uses \texttt{bfindzero} and \texttt{ratezero} to find the 
zero entry most likely to complete a row, column, or Reidemeister III word, 
then systematically fills in the
zero with each possible entry returned by \texttt{avail}. For each such value,
the pattern is run through \texttt{biqfill} and any patterns returned are 
added to the end of the working list. When a pattern has no remaining 
zeroes, it is transferred to the output list; when the working list is empty, 
the program reports the output list. \texttt{biqlist} is useful for finding 
biquandle completions for knot biquandles, and we used it to find all 
biquandles of order 2, 3 and 4 as well as several biquandles of orders 5 and 6.

To compute the biquandle counting invariant, we use \texttt{bhomlist} which
represents a map $\phi:\{1,2,\dots, n\}\to \{1,2,\dots,m\}$ as
a vector $[\phi(1),\phi(2),\dots,\phi(n)]$ where $\phi(i)\in \{1,2,\dots, m\}.$
This program starts with a vector of all zeroes and keeps a working list of 
vectors, systematically filling in zeroes with each possible entry and using 
\texttt{bhomfill} to propagate the 
values through the vector as required by the homomorphism conditions, removing
any maps which violate the homomorphism conditions. Once all zeroes have been 
filled in, each completed homomorphism is moved to the output list, and once
the working list is empty, the output list is reported. The program 
\texttt{bhomcount} reports the cardinality of the list output by 
\texttt{bhomlist}, that is, the value of the counting invariant for the
input presentation and target biquandles. 

\texttt{bhomlist} can compute the set of homomorphisms from any biquandle 
presentation matrix to any finite biquandle. Since a finite biquandle matrix 
is itself a presentation matrix, we can use \texttt{bhomlist} to determine 
whether any pair of finite biquandles are isomorphic and to compute the 
automorphism group of any finite biquandle. \texttt{bisolist} lists the
isomorphisms from one finite biquandle to another, giving an empty list if
the two are non-isomorphic. \texttt{baut} computes the automorphism group
of a finite biquandle.

\texttt{Obv} gives the obverse of a biquandle matrix, \texttt{Flip} gives the
flip. The program \texttt{breducelist} takes a list of biquandles and
removes any which are isomorphic to previous biquandles on the list, their 
flips, their obverses, or their obverse flips. For long lists of biquandles,
it may be necessary to split the list into smaller lists first.

We used \texttt{biqlist} and \texttt{breducelist} to classify finite 
biquandles of order 2, 3 and 4; the results are lists in tables \ref{t1}
through \ref{t8}. There are 36 biquandles of order 3, which comprise 15 
isomorphism classes. Five of these
are isomorphic to the flip, obverse, or obverse flip of one of the ten
listed biquandles of order 3. There are 744 biquandles of order 4, which 
reduce to 64 when isomorphic copies, flips and obverses are removed. In the 
interest of space, we have listed only the non-quandle biquandles of order 4.
In the tables we have noted whether a biquandle is self-flip and listed the 
automorphism groups computed with \texttt{baut}.

\begin{table} 
\begin{center}
\begin{tabular}{|ccc|ccc|}  \hline
  & & & & &  \\
Biquandle  & Self- & Automorphism & Biquandle  & Self- & Automorphism \\
Matrix  & Flip? & Group & Matrix  & Flip? & Group \\
 & & & & &  \\ \hline
 & & & & &  \\
$\displaystyle{\left[
\begin{array}{cc|cc}
1 & 1 & 1 & 1 \\
2 & 2 & 2 & 2 \\ \hline
1 & 1 & 1 & 1 \\
2 & 2 & 2 & 2 \\
\end{array}\right]}$  & Yes & $\mathbb{Z}_2$ &
$\displaystyle{\left[
\begin{array}{cc|cc}
2 & 2 & 2 & 2 \\ 
1 & 1 & 1 & 1 \\ \hline
2 & 2 & 2 & 2 \\
1 & 1 & 1 & 1 \\
\end{array}\right]}$  & Yes & $\mathbb{Z}_2$ \\
 & & & & & \\
\hline
\end{tabular}
\end{center}
\caption{Biquandles of order 2.} 
\label{t1}
\end{table}

\begin{table}
\begin{center}
\begin{tabular}{|ccc|ccc|}  \hline
  & & & & & \\
Biquandle & Self- & $\mathrm{Aut}(B)$ & Biquandle & Self- & $\mathrm{Aut}(B)$\\
Matrix  & Flip? &  & Matrix  & Flip? &  \\
 & & & & & \\ \hline
 & & & & & \\
$\displaystyle{\left[
\begin{array}{ccc|ccc}
1 & 1 & 1 & 1 & 1 & 1 \\
3 & 3 & 3 & 3 & 3 & 3 \\
2 & 2 & 2 & 2 & 2 & 2 \\ \hline
1 & 3 & 2 & 1 & 3 & 2 \\
2 & 1 & 3 & 2 & 1 & 3 \\
3 & 2 & 1 & 3 & 2 & 1
\end{array}
\right]}$ & No & $\mathbb{Z}_2$  &
$\displaystyle{\left[
\begin{array}{ccc|ccc}
1 & 3 & 2 & 1 & 3 & 2 \\
3 & 2 & 1 & 3 & 2 & 1 \\
2 & 1 & 3 & 2 & 1 & 3 \\ 
\hline
1 & 1 & 1 & 1 & 1 & 1 \\
2 & 2 & 2 & 2 & 2 & 2 \\
3 & 3 & 3 & 3 & 3 & 3
\end{array}
\right]}$  & No & $S_3$  \\
 & & & & & \\
$\displaystyle{\left[
\begin{array}{ccc|ccc}
2 & 1 & 3 & 3 & 2 & 1 \\
1 & 3 & 2 & 2 & 1 & 3 \\
3 & 2 & 1 & 1 & 3 & 2 \\
\hline
2 & 2 & 2 & 3 & 3 & 3 \\
3 & 3 & 3 & 1 & 1 & 1 \\
1 & 1 & 1 & 2 & 2 & 2 \\ 
\end{array}
\right]}$ & No & $\mathbb{Z}_3$ &
$\displaystyle{\left[
\begin{array}{ccc|ccc}
2 & 2 & 2 & 3 & 3 & 3 \\
3 & 3 & 3 & 1 & 1 & 1 \\
1 & 1 & 1 & 2 & 2 & 2 \\ 
\hline
3 & 3 & 3 & 2 & 2 & 2 \\
1 & 1 & 1 & 3 & 3 & 3 \\
2 & 2 & 2 & 1 & 1 & 1 \\
\end{array}
\right]}$ & Yes & $\mathbb{Z}_3$ \\
 & & & & & \\
$\displaystyle{\left[
\begin{array}{ccc|ccc}
1 & 1 & 1 & 1 & 1 & 1 \\
2 & 2 & 2 & 2 & 2 & 2 \\
3 & 3 & 3 & 3 & 3 & 3 \\ 
\hline
1 & 1 & 2 & 1 & 1 & 2 \\
2 & 2 & 1 & 2 & 2 & 1 \\
3 & 3 & 3 & 3 & 3 & 3 \\ 
\end{array}
\right]}$  & No & $\mathbb{Z}_2$ &
$\displaystyle{\left[
\begin{array}{ccc|ccc}
1 & 1 & 1 & 1 & 1 & 1 \\
2 & 3 & 3 & 2 & 3 & 3 \\
3 & 2 & 2 & 3 & 2 & 2 \\ 
\hline
1 & 1 & 1 & 1 & 1 & 1 \\
2 & 3 & 3 & 2 & 3 & 3 \\
3 & 2 & 2 & 3 & 2 & 2 \\
\end{array}
\right]}$  & Yes & $\mathbb{Z}_2$ \\
& & & & & \\
$\displaystyle{\left[
\begin{array}{ccc|ccc}
1 & 1 & 1 & 1 & 1 & 1 \\
2 & 3 & 3 & 2 & 3 & 3 \\
3 & 2 & 2 & 3 & 2 & 2 \\ 
\hline
1 & 1 & 1 & 1 & 1 & 1 \\
3 & 3 & 3 & 3 & 3 & 3 \\
2 & 2 & 2 & 2 & 2 & 2 \\
\end{array}
\right]}$  & No & $\mathbb{Z}_2$  & 
$\displaystyle{\left[
\begin{array}{ccc|ccc}
1 & 1 & 1 & 1 & 1 & 1 \\
3 & 2 & 2 & 3 & 2 & 2 \\
2 & 3 & 3 & 2 & 3 & 3 \\ 
\hline
1 & 1 & 1 & 1 & 1 & 1 \\
3 & 2 & 2 & 3 & 2 & 2 \\
2 & 3 & 3 & 2 & 3 & 3 \\
\end{array}
\right]}$  & Yes & $\mathbb{Z}_2$ \\
& & & & & \\
$\displaystyle{\left[
\begin{array}{ccc|ccc}
1 & 1 & 1 & 1 & 1 & 1 \\
3 & 3 & 3 & 3 & 3 & 3 \\
2 & 2 & 2 & 2 & 2 & 2 \\ 
\hline
1 & 1 & 1 & 1 & 1 & 1 \\
3 & 3 & 3 & 3 & 3 & 3 \\
2 & 2 & 2 & 2 & 2 & 2 \\
\end{array}
\right]}$ & Yes & $\mathbb{Z}_2$ &
$\displaystyle{\left[
\begin{array}{ccc|ccc}
1 & 1 & 1 & 1 & 1 & 1 \\
2 & 2 & 2 & 2 & 2 & 2 \\ 
3 & 3 & 3 & 3 & 3 & 3 \\
\hline
1 & 1 & 1 & 1 & 1 & 1 \\
2 & 2 & 2 & 2 & 2 & 2 \\
3 & 3 & 3 & 3 & 3 & 3 \\
\end{array}
\right]}$  & Yes & $S_3$ \\
& & & & & \\
\hline
\end{tabular}
\end{center}
\caption{Biquandles of order 3.}
\label{t2}
\end{table}

\begin{table} \footnotesize{
\begin{center}
\begin{tabular}{|ccc|ccc|}  \hline
  & & & & & \\
Biquandle  & Self- & $\mathrm{Aut}(B)$ &
Biquandle  & Self- & $\mathrm{Aut}(B)$\\
Matrix  & Flip?  & & Matrix & Flip? & \\
 & & & & & \\ \hline
 & & & & & \\
$\displaystyle{\left[
\begin{array}{cccc|cccc}
2 & 2 & 2 & 2 & 2 & 2 & 2 & 2 \\
1 & 1 & 1 & 1 & 1 & 1 & 1 & 1 \\
4 & 4 & 4 & 4 & 4 & 4 & 4 & 4 \\
3 & 3 & 3 & 3 & 3 & 3 & 3 & 3 \\ \hline
3 & 2 & 2 & 3 & 3 & 2 & 2 & 3 \\
1 & 4 & 4 & 1 & 1 & 4 & 4 & 1 \\
4 & 1 & 1 & 4 & 4 & 1 & 1 & 4 \\
2 & 3 & 3 & 2 & 2 & 3 & 3 & 2 \\
\end{array}\right]}$ & No & $\mathbb{Z}_2\oplus\mathbb{Z}_2$ &
$\displaystyle{\left[
\begin{array}{cccc|cccc}
2 & 4 & 3 & 1 & 4 & 1 & 3 & 2 \\
3 & 1 & 2 & 4 & 2 & 3 & 1 & 4 \\
1 & 3 & 4 & 2 & 1 & 4 & 2 & 3 \\
4 & 2 & 1 & 3 & 3 & 2 & 4 & 1 \\ \hline
3 & 3 & 3 & 3 & 3 & 3 & 3 & 3 \\
4 & 4 & 4 & 4 & 4 & 4 & 4 & 4 \\
1 & 1 & 1 & 1 & 1 & 1 & 1 & 1 \\
2 & 2 & 2 & 2 & 2 & 2 & 2 & 2 \\
\end{array}\right]}$ & No & $\mathbb{Z}_2\oplus\mathbb{Z}_2$ \\
& & & & & \\
$\displaystyle{\left[
\begin{array}{cccc|cccc}
1 & 4 & 2 & 3 & 1 & 3 & 4 & 2 \\
2 & 3 & 1 & 4 & 3 & 1 & 2 & 4 \\
4 & 1 & 3 & 2 & 2 & 4 & 3 & 1 \\
3 & 2 & 4 & 1 & 4 & 2 & 1 & 3 \\ \hline
1 & 3 & 4 & 2 & 1 & 4 & 2 & 3 \\
3 & 1 & 2 & 4 & 2 & 3 & 1 & 4 \\
2 & 4 & 3 & 1 & 4 & 1 & 3 & 2 \\
4 & 2 & 1 & 3 & 3 & 2 & 4 & 1 \\ 
\end{array}\right]}$ & Yes & $\mathbb{Z}_2$ &
$\displaystyle{\left[
\begin{array}{cccc|cccc}
1 & 4 & 2 & 3 & 1 & 3 & 4 & 2 \\
2 & 3 & 1 & 4 & 3 & 1 & 2 & 4 \\
4 & 1 & 3 & 2 & 2 & 4 & 3 & 1 \\
3 & 2 & 4 & 1 & 4 & 2 & 1 & 3 \\ \hline
1 & 1 & 1 & 1 & 1 & 1 & 1 & 1 \\
3 & 3 & 3 & 3 & 4 & 4 & 4 & 4 \\
4 & 4 & 4 & 4 & 2 & 2 & 2 & 2 \\
2 & 2 & 2 & 2 & 3 & 3 & 3 & 3 \\ 
\end{array}\right]}$ & No & $\mathbb{Z}_3$ \\
& & & & & \\
$\displaystyle{\left[
\begin{array}{cccc|cccc}
1 & 3 & 4 & 2 & 1 & 4 & 2 & 3 \\
3 & 1 & 2 & 4 & 4 & 1 & 3 & 2 \\
4 & 2 & 1 & 3 & 2 & 3 & 1 & 4 \\
2 & 4 & 3 & 1 & 3 & 2 & 4 & 1 \\ \hline
1 & 1 & 1 & 1 & 1 & 1 & 1 & 1 \\
3 & 3 & 3 & 3 & 4 & 4 & 4 & 4 \\
4 & 4 & 4 & 4 & 2 & 2 & 2 & 2 \\
2 & 2 & 2 & 2 & 3 & 3 & 3 & 3 \\ 
\end{array}\right]}$ & No & $\mathbb{Z}_3$ &
$\displaystyle{\left[
\begin{array}{cccc|cccc}
3 & 2 & 1 & 2 & 3 & 2 & 1 & 2 \\
1 & 3 & 2 & 1 & 1 & 3 & 2 & 1 \\
2 & 1 & 3 & 3 & 2 & 1 & 3 & 3 \\
4 & 4 & 4 & 4 & 4 & 4 & 4 & 4 \\ \hline
2 & 2 & 2 & 2 & 2 & 2 & 2 & 2 \\ 
1 & 1 & 1 & 1 & 1 & 1 & 1 & 1 \\
3 & 3 & 3 & 3 & 3 & 3 & 3 & 3 \\
4 & 4 & 4 & 4 & 4 & 4 & 4 & 4 \\
\end{array}\right]}$ & No & $\mathbb{Z}_2$ \\
& & & & & \\
$\displaystyle{\left[
\begin{array}{cccc|cccc}
3 & 2 & 1 & 1 & 3 & 2 & 1 & 1 \\
1 & 3 & 2 & 2 & 1 & 3 & 2 & 2 \\
2 & 1 & 3 & 3 & 2 & 1 & 3 & 3 \\
4 & 4 & 4 & 4 & 4 & 4 & 4 & 4 \\ \hline
2 & 2 & 2 & 1 & 2 & 2 & 2 & 1 \\ 
1 & 1 & 1 & 2 & 1 & 1 & 1 & 2 \\
3 & 3 & 3 & 3 & 3 & 3 & 3 & 3 \\
4 & 4 & 4 & 4 & 4 & 4 & 4 & 4 \\
\end{array}\right]}$ & No & $\mathbb{Z}_3$ &
$\displaystyle{\left[
\begin{array}{cccc|cccc}
3 & 2 & 4 & 1 & 3 & 2 & 1 & 4 \\
4 & 1 & 3 & 2 & 1 & 3 & 4 & 2 \\
2 & 3 & 1 & 4 & 4 & 2 & 1 & 3 \\
1 & 4 & 2 & 3 & 2 & 4 & 3 & 1 \\ \hline
3 & 1 & 2 & 4 & 3 & 2 & 4 & 1 \\ 
1 & 3 & 4 & 2 & 4 & 1 & 3 & 2 \\
4 & 2 & 1 & 3 & 2 & 3 & 1 & 4 \\
2 & 4 & 3 & 1 & 1 & 4 & 2 & 3 \\
\end{array}\right]}$ & Yes & $\mathbb{Z}_2$ \\
& & & & & \\
$\displaystyle{\left[
\begin{array}{cccc|cccc}
1 & 1 & 1 & 2 & 1 & 1 & 1 & 3 \\
2 & 2 & 2 & 3 & 2 & 2 & 2 & 1 \\
3 & 3 & 3 & 1 & 3 & 3 & 3 & 2 \\
4 & 4 & 4 & 4 & 4 & 4 & 4 & 4 \\ \hline
1 & 3 & 2 & 3 & 1 & 3 & 2 & 2 \\ 
3 & 2 & 1 & 1 & 3 & 2 & 1 & 3 \\
2 & 1 & 3 & 2 & 2 & 1 & 3 & 1 \\
4 & 4 & 4 & 4 & 4 & 4 & 4 & 4 \\
\end{array}\right]}$ & No & $\mathbb{Z}_3$ &
$\displaystyle{\left[
\begin{array}{cccc|cccc}
1 & 1 & 1 & 2 & 1 & 1 & 1 & 2 \\
2 & 2 & 2 & 3 & 2 & 2 & 2 & 1 \\
3 & 3 & 3 & 1 & 3 & 3 & 3 & 3 \\
4 & 4 & 4 & 4 & 4 & 4 & 4 & 4 \\ \hline
1 & 3 & 2 & 2 & 1 & 3 & 2 & 2 \\ 
3 & 2 & 1 & 1 & 3 & 2 & 1 & 1 \\
2 & 1 & 3 & 3 & 2 & 1 & 3 & 3 \\
4 & 4 & 4 & 4 & 4 & 4 & 4 & 4 \\
\end{array}\right]}$ & No & $\mathbb{Z}_2$ \\
& & & & & \\
\hline
\end{tabular}
\end{center}}
\caption{Non-quandle Biquandles of order 4 part 1.} 
\label{t3}
\end{table}

\begin{table} \footnotesize{
\begin{center}
\begin{tabular}{|ccc|ccc|}  \hline
  & & & & & \\
Biquandle  & Self- & $\mathrm{Aut}(B)$ &
Biquandle & Self- & $\mathrm{Aut}(B)$\\
Matrix  & Flip?  & & Matrix & Flip? & \\
 & & & & & \\ \hline
 & & & & & \\
$\displaystyle{\left[
\begin{array}{cccc|cccc}
2 & 2 & 2 & 3 & 3 & 3 & 3 & 2 \\
3 & 3 & 3 & 1 & 1 & 1 & 3 & 2 \\
1 & 1 & 1 & 2 & 2 & 2 & 2 & 1 \\
4 & 4 & 4 & 4 & 4 & 4 & 4 & 4 \\ \hline
3 & 3 & 3 & 3 & 2 & 2 & 2 & 2 \\
1 & 1 & 1 & 1 & 3 & 3 & 3 & 3 \\
2 & 2 & 2 & 2 & 1 & 1 & 1 & 1 \\
4 & 4 & 4 & 4 & 4 & 4 & 4 & 4 \\
\end{array}\right]}$ & No & $\mathbb{Z}_3$ &
$\displaystyle{\left[
\begin{array}{cccc|cccc}
2 & 2 & 2 & 3 & 3 & 3 & 3 & 2 \\
3 & 3 & 3 & 1 & 1 & 1 & 1 & 3 \\
1 & 1 & 1 & 2 & 2 & 2 & 2 & 1 \\
4 & 4 & 4 & 4 & 4 & 4 & 4 & 4 \\ \hline
3 & 3 & 3 & 1 & 2 & 2 & 2 & 1 \\
1 & 1 & 1 & 2 & 3 & 3 & 3 & 2 \\
2 & 2 & 2 & 3 & 1 & 1 & 1 & 3 \\
4 & 4 & 4 & 4 & 4 & 4 & 4 & 4 \\
\end{array}\right]}$ & No & $\mathbb{Z}_3$ \\
& & & & & \\
$\displaystyle{\left[
\begin{array}{cccc|cccc}
2 & 2 & 2 & 3 & 3 & 3 & 3 & 2 \\
3 & 3 & 3 & 1 & 1 & 1 & 1 & 3 \\
1 & 1 & 1 & 2 & 2 & 2 & 2 & 1 \\
4 & 4 & 4 & 4 & 4 & 4 & 4 & 4 \\ \hline
3 & 3 & 3 & 2 & 2 & 2 & 2 & 3 \\
1 & 1 & 1 & 3 & 3 & 3 & 3 & 1 \\
2 & 2 & 2 & 1 & 1 & 1 & 1 & 2 \\
4 & 4 & 4 & 4 & 4 & 4 & 4 & 4 \\
\end{array}\right]}$ & Yes & $\mathbb{Z}_3$ &
$\displaystyle{\left[
\begin{array}{cccc|cccc}
2 & 2 & 2 & 1 & 3 & 3 & 3 & 1 \\
3 & 3 & 3 & 2 & 1 & 1 & 1 & 2 \\
1 & 1 & 1 & 3 & 2 & 2 & 2 & 3 \\
4 & 4 & 4 & 4 & 4 & 4 & 4 & 4 \\ \hline
3 & 3 & 3 & 3 & 2 & 2 & 2 & 2 \\
1 & 1 & 1 & 1 & 3 & 3 & 3 & 3 \\
2 & 2 & 2 & 2 & 1 & 1 & 1 & 1 \\
4 & 4 & 4 & 4 & 4 & 4 & 4 & 4 \\
\end{array}\right]}$ & No & $\mathbb{Z}_3$ \\
& & & & & \\
$\displaystyle{\left[
\begin{array}{cccc|cccc}
2 & 2 & 2 & 1 & 3 & 3 & 3 & 1 \\
3 & 3 & 3 & 2 & 1 & 1 & 1 & 2 \\
1 & 1 & 1 & 3 & 2 & 2 & 2 & 3 \\
4 & 4 & 4 & 4 & 4 & 4 & 4 & 4 \\ \hline
3 & 3 & 3 & 1 & 2 & 2 & 2 & 1 \\
1 & 1 & 1 & 2 & 3 & 3 & 3 & 2 \\
2 & 2 & 2 & 3 & 1 & 1 & 1 & 3 \\
4 & 4 & 4 & 4 & 4 & 4 & 4 & 4 \\
\end{array}\right]}$ & Yes & $\mathbb{Z}_3$ &
$\displaystyle{\left[
\begin{array}{cccc|cccc}
2 & 2 & 2 & 2 & 3 & 3 & 3 & 3 \\
3 & 3 & 3 & 3 & 1 & 1 & 1 & 1 \\
1 & 1 & 1 & 1 & 2 & 2 & 2 & 2 \\
4 & 4 & 4 & 4 & 4 & 4 & 4 & 4 \\ \hline
3 & 3 & 3 & 3 & 2 & 2 & 2 & 2 \\
1 & 1 & 1 & 1 & 3 & 3 & 3 & 3 \\
2 & 2 & 2 & 2 & 1 & 1 & 1 & 1 \\
4 & 4 & 4 & 4 & 4 & 4 & 4 & 4 \\
\end{array}\right]}$ & Yes & $\mathbb{Z}_3$ \\
& & & & & \\
$\displaystyle{\left[
\begin{array}{cccc|cccc}
3 & 2 & 2 & 3 & 3 & 2 & 2 & 3 \\
1 & 4 & 4 & 1 & 1 & 4 & 4 & 1 \\
4 & 1 & 1 & 4 & 4 & 1 & 1 & 4 \\
2 & 3 & 3 & 2 & 2 & 3 & 3 & 2 \\ \hline
3 & 2 & 2 & 3 & 3 & 2 & 2 & 3 \\
1 & 4 & 4 & 1 & 1 & 4 & 4 & 1 \\
4 & 1 & 1 & 4 & 4 & 1 & 1 & 4 \\
2 & 3 & 3 & 2 & 2 & 3 & 3 & 2 \\
\end{array}\right]}$ & Yes & $\mathbb{Z}_2\oplus \mathbb{Z}_2$ &
$\displaystyle{\left[
\begin{array}{cccc|cccc}
3 & 2 & 1 & 3 & 2 & 1 & 3 & 2 \\
2 & 1 & 3 & 1 & 1 & 3 & 2 & 3 \\
1 & 3 & 2 & 2 & 3 & 2 & 1 & 1 \\
4 & 4 & 4 & 4 & 4 & 4 & 4 & 4 \\ \hline
3 & 3 & 3 & 2 & 2 & 2 & 2 & 3 \\
1 & 1 & 1 & 3 & 3 & 3 & 3 & 1 \\
2 & 2 & 2 & 1 & 1 & 1 & 1 & 2 \\
4 & 4 & 4 & 4 & 4 & 4 & 4 & 4 \\
\end{array}\right]}$ & No & $\mathbb{Z}_3$ \\
& & & & & \\
$\displaystyle{\left[
\begin{array}{cccc|cccc}
3 & 2 & 1 & 1 & 2 & 1 & 2 & 1 \\
2 & 1 & 3 & 2 & 1 & 3 & 2 & 2 \\
1 & 3 & 2 & 3 & 3 & 2 & 1 & 3 \\
4 & 4 & 4 & 4 & 4 & 4 & 4 & 4 \\ \hline
3 & 3 & 3 & 1 & 2 & 2 & 2 & 1 \\
1 & 1 & 1 & 2 & 3 & 3 & 3 & 2 \\
2 & 2 & 2 & 3 & 1 & 1 & 1 & 3 \\
4 & 4 & 4 & 4 & 4 & 4 & 4 & 4 \\
\end{array}\right]}$ & No & $\mathbb{Z}_3$ &
$\displaystyle{\left[
\begin{array}{cccc|cccc}
3 & 2 & 1 & 2 & 2 & 1 & 3 & 3 \\
2 & 1 & 3 & 3 & 1 & 3 & 2 & 1 \\
1 & 3 & 2 & 1 & 3 & 2 & 1 & 3 \\
4 & 4 & 4 & 4 & 4 & 4 & 4 & 4 \\ \hline
3 & 3 & 3 & 3 & 2 & 2 & 2 & 2 \\
1 & 1 & 1 & 1 & 3 & 3 & 3 & 3 \\
2 & 2 & 2 & 2 & 1 & 1 & 1 & 1 \\
4 & 4 & 4 & 4 & 4 & 4 & 4 & 4 \\
\end{array}\right]}$ & No & $\mathbb{Z}_3$ \\
& & & & & \\
\hline
\end{tabular}
\end{center}}
\caption{Non-quandle Biquandles of order 4 part 2.} 
\label{t4}
\end{table}

\begin{table} \footnotesize{
\begin{center}
\begin{tabular}{|ccc|ccc|}  \hline
  & & & & & \\
Biquandle & Self- & $\mathrm{Aut}(B)$ &
Biquandle & Self- & $\mathrm{Aut}(B)$\\
Matrix  & Flip?  & & Matrix & Flip? & \\
 & & & & & \\ \hline
 & & & & & \\
$\displaystyle{\left[
\begin{array}{cccc|cccc}
2 & 3 & 3 & 2 & 3 & 2 & 2 & 3 \\
1 & 4 & 4 & 1 & 4 & 1 & 1 & 4 \\
4 & 1 & 1 & 4 & 1 & 4 & 4 & 1 \\
3 & 2 & 2 & 3 & 2 & 3 & 3 & 2 \\ \hline
2 & 2 & 2 & 2 & 3 & 3 & 3 & 3 \\
4 & 4 & 4 & 4 & 1 & 1 & 1 & 1 \\
1 & 1 & 1 & 1 & 4 & 4 & 4 & 4 \\
3 & 3 & 3 & 3 & 2 & 2 & 2 & 2 \\
\end{array}\right]}$ & No & $\mathbb{Z}_4$ &
$\displaystyle{\left[
\begin{array}{cccc|cccc}
4 & 4 & 3 & 3 & 3 & 3 & 4 & 4 \\
3 & 3 & 4 & 4 & 4 & 4 & 3 & 3 \\
2 & 2 & 1 & 1 & 1 & 1 & 2 & 2 \\
1 & 1 & 2 & 2 & 2 & 2 & 1 & 1 \\ \hline
3 & 3 & 4 & 4 & 4 & 4 & 3 & 3 \\
4 & 4 & 3 & 3 & 3 & 3 & 4 & 4 \\
1 & 1 & 2 & 2 & 2 & 2 & 1 & 1 \\
2 & 2 & 1 & 1 & 1 & 1 & 2 & 2 \\
\end{array}\right]}$ & Yes & $\mathbb{Z}_4 $ \\
& & & & & \\
$\displaystyle{\left[
\begin{array}{cccc|cccc}
4 & 4 & 4 & 4 & 3 & 3 & 3 & 3 \\
3 & 3 & 3 & 3 & 4 & 4 & 4 & 4 \\
1 & 1 & 1 & 1 & 2 & 2 & 2 & 2 \\
2 & 2 & 2 & 2 & 1 & 1 & 1 & 1 \\ \hline
3 & 3 & 3 & 3 & 4 & 4 & 4 & 4 \\
4 & 4 & 4 & 4 & 3 & 3 & 3 & 3 \\
2 & 2 & 2 & 2 & 1 & 1 & 1 & 1 \\
1 & 1 & 1 & 1 & 2 & 2 & 2 & 2 \\
\end{array}\right]}$ & Yes & $\mathbb{Z}_4$ &
$\displaystyle{\left[
\begin{array}{cccc|cccc}
1 & 1 & 1 & 2 & 1 & 1 & 1 & 3 \\
2 & 2 & 2 & 3 & 2 & 2 & 2 & 1 \\
3 & 3 & 3 & 1 & 3 & 3 & 3 & 2 \\
4 & 4 & 4 & 4 & 4 & 4 & 4 & 4 \\ \hline
1 & 1 & 1 & 3 & 1 & 1 & 1 & 2 \\
2 & 2 & 2 & 1 & 2 & 2 & 2 & 3 \\
3 & 3 & 3 & 2 & 3 & 3 & 3 & 1 \\
4 & 4 & 4 & 4 & 4 & 4 & 4 & 4 \\
\end{array}\right]}$ & Yes & $\mathbb{Z}_3$ \\
& & & & & \\
$\displaystyle{\left[
\begin{array}{cccc|cccc}
2 & 2 & 1 & 1 & 2 & 2 & 1 & 1 \\
1 & 1 & 2 & 2 & 1 & 1 & 2 & 2 \\
4 & 3 & 3 & 3 & 4 & 3 & 3 & 3 \\
3 & 4 & 4 & 4 & 3 & 4 & 4 & 4 \\ \hline
2 & 2 & 1 & 1 & 2 & 2 & 1 & 1 \\
1 & 1 & 2 & 2 & 1 & 1 & 2 & 2 \\
3 & 4 & 3 & 3 & 3 & 4 & 3 & 3 \\
4 & 3 & 4 & 4 & 4 & 3 & 4 & 4 \\
\end{array}\right]}$ & Yes & $\mathbb{Z}_2$ &
$\displaystyle{\left[
\begin{array}{cccc|cccc}
2 & 2 & 1 & 1 & 2 & 2 & 1 & 1 \\
1 & 1 & 2 & 2 & 1 & 1 & 2 & 2 \\
4 & 3 & 4 & 4 & 4 & 3 & 4 & 4 \\
3 & 4 & 3 & 3 & 3 & 4 & 3 & 3 \\ \hline
2 & 2 & 1 & 1 & 2 & 2 & 1 & 1 \\
1 & 1 & 2 & 2 & 1 & 1 & 2 & 2 \\
3 & 4 & 4 & 4 & 3 & 4 & 4 & 4 \\
4 & 3 & 3 & 3 & 4 & 3 & 3 & 3 \\
\end{array}\right]}$ & Yes & $\mathbb{Z}_2$ \\
& & & & & \\
$\displaystyle{\left[
\begin{array}{cccc|cccc}
2 & 2 & 1 & 1 & 2 & 2 & 1 & 1 \\
1 & 1 & 2 & 2 & 1 & 1 & 2 & 2 \\
3 & 4 & 3 & 3 & 3 & 4 & 3 & 3 \\
4 & 3 & 4 & 4 & 4 & 3 & 4 & 4 \\ \hline
2 & 2 & 1 & 1 & 2 & 2 & 1 & 1 \\
1 & 1 & 2 & 2 & 1 & 1 & 2 & 2 \\
3 & 4 & 3 & 3 & 3 & 4 & 3 & 3 \\
4 & 3 & 4 & 4 & 4 & 3 & 4 & 4 \\
\end{array}\right]}$ & Yes & $\mathbb{Z}_2$ &
$\displaystyle{\left[
\begin{array}{cccc|cccc}
2 & 2 & 1 & 1 & 2 & 2 & 1 & 1 \\
1 & 1 & 2 & 2 & 1 & 1 & 2 & 1 \\
3 & 4 & 4 & 4 & 3 & 4 & 4 & 4 \\
4 & 3 & 3 & 3 & 4 & 3 & 3 & 3 \\ \hline
2 & 2 & 1 & 1 & 2 & 2 & 1 & 1 \\
1 & 1 & 2 & 2 & 1 & 1 & 2 & 2 \\
3 & 4 & 4 & 4 & 3 & 4 & 4 & 4 \\
4 & 3 & 3 & 3 & 4 & 3 & 3 & 3 \\
\end{array}\right]}$ & Yes & $\mathbb{Z}_2$ \\
& & & & & \\
$\displaystyle{\left[
\begin{array}{cccc|cccc}
1 & 1 & 1 & 1 & 1 & 1 & 1 & 1 \\
3 & 2 & 2 & 3 & 3 & 2 & 2 & 3 \\
2 & 3 & 3 & 2 & 2 & 3 & 3 & 2 \\
4 & 4 & 4 & 4 & 4 & 4 & 4 & 4 \\ \hline
1 & 1 & 1 & 1 & 1 & 1 & 1 & 1 \\
2 & 2 & 2 & 3 & 2 & 2 & 2 & 3 \\
3 & 3 & 3 & 2 & 3 & 3 & 3 & 2 \\
4 & 4 & 4 & 4 & 4 & 4 & 4 & 4 \\
\end{array}\right]}$ & No & $\mathbb{Z}_2$ &
$\displaystyle{\left[
\begin{array}{cccc|cccc}
1 & 1 & 1 & 1 & 1 & 1 & 1 & 1 \\
3 & 2 & 2 & 2 & 4 & 2 & 2 & 2 \\
4 & 3 & 3 & 3 & 2 & 3 & 3 & 3 \\
2 & 4 & 4 & 4 & 3 & 4 & 4 & 4 \\ \hline
1 & 1 & 1 & 1 & 1 & 1 & 1 & 1 \\
3 & 2 & 2 & 2 & 4 & 2 & 2 & 2 \\
4 & 3 & 3 & 3 & 2 & 3 & 3 & 3 \\
2 & 4 & 4 & 4 & 3 & 4 & 4 & 4 \\
\end{array}\right]}$ & Yes & $\mathbb{Z}_3$ \\
& & & & & \\
\hline
\end{tabular}
\end{center}
\caption{Non-quandle Biquandles of order 4 part 3.} }
\label{t5}
\end{table}

\begin{table} \footnotesize{
\begin{center}
\begin{tabular}{|ccc|ccc|}  \hline
  & & & & & \\
Biquandle & Self- & $\mathrm{Aut}(B)$ &
Biquandle & Self- & $\mathrm{Aut}(B)$\\
Matrix  & Flip?  & & Matrix & Flip? & \\
 & & & & & \\ \hline
 & & & & & \\
$\displaystyle{\left[
\begin{array}{cccc|cccc}
1 & 1 & 1 & 1 & 1 & 1 & 1 & 1 \\
3 & 2 & 2 & 2 & 3 & 2 & 2 & 2 \\
2 & 3 & 3 & 3 & 2 & 3 & 3 & 3 \\
4 & 4 & 4 & 4 & 4 & 4 & 4 & 4 \\ \hline
1 & 1 & 1 & 1 & 1 & 1 & 1 & 1 \\
2 & 2 & 2 & 3 & 2 & 2 & 2 & 3 \\
3 & 3 & 3 & 2 & 3 & 3 & 3 & 2 \\
4 & 4 & 4 & 4 & 4 & 4 & 4 & 4 \\
\end{array}\right]}$ & Yes & $\mathbb{Z}_2$ &
$\displaystyle{\left[
\begin{array}{cccc|cccc}
1 & 1 & 1 & 1 & 1 & 1 & 1 & 1 \\
3 & 3 & 3 & 3 & 3 & 3 & 3 & 3 \\
2 & 2 & 2 & 2 & 2 & 2 & 2 & 2 \\
4 & 4 & 4 & 4 & 4 & 4 & 4 & 4 \\ \hline
1 & 1 & 1 & 1 & 1 & 1 & 1 & 1 \\
2 & 3 & 3 & 3 & 2 & 3 & 3 & 3 \\
3 & 2 & 2 & 2 & 3 & 2 & 2 & 2 \\
4 & 4 & 4 & 4 & 4 & 4 & 4 & 4 \\
\end{array}\right]}$ & No & $\mathbb{Z}_2$ \\
& & & & & \\
$\displaystyle{\left[
\begin{array}{cccc|cccc}
1 & 1 & 1 & 1 & 1 & 1 & 1 & 1 \\
3 & 3 & 3 & 3 & 3 & 3 & 3 & 3 \\
2 & 2 & 2 & 2 & 2 & 2 & 2 & 2 \\
4 & 4 & 4 & 4 & 4 & 4 & 4 & 4 \\ \hline
1 & 1 & 1 & 1 & 1 & 1 & 1 & 1 \\
2 & 3 & 3 & 2 & 2 & 3 & 3 & 2 \\
3 & 2 & 2 & 3 & 3 & 2 & 2 & 3 \\
4 & 4 & 4 & 4 & 4 & 4 & 4 & 4 \\
\end{array}\right]}$ & No & $\mathbb{Z}_2\oplus \mathbb{Z}_2$ &
$\displaystyle{\left[
\begin{array}{cccc|cccc}
1 & 1 & 1 & 1 & 1 & 1 & 1 & 1 \\
3 & 3 & 3 & 2 & 3 & 3 & 3 & 2 \\
2 & 2 & 2 & 3 & 2 & 2 & 2 & 3 \\
4 & 4 & 4 & 4 & 4 & 4 & 4 & 4 \\ \hline
1 & 1 & 1 & 1 & 1 & 1 & 1 & 1 \\
2 & 3 & 3 & 3 & 2 & 3 & 3 & 3 \\
3 & 2 & 2 & 2 & 3 & 2 & 2 & 2 \\
4 & 4 & 4 & 4 & 4 & 4 & 4 & 4 \\
\end{array}\right]}$ & Yes & $\mathbb{Z}_2$ \\
& & & & & \\
$\displaystyle{\left[
\begin{array}{cccc|cccc}
1 & 1 & 1 & 1 & 1 & 1 & 1 & 1 \\
3 & 3 & 3 & 2 & 3 & 3 & 3 & 2 \\
2 & 2 & 2 & 3 & 2 & 2 & 2 & 3 \\
4 & 4 & 4 & 4 & 4 & 4 & 4 & 4 \\ \hline
1 & 1 & 1 & 1 & 1 & 1 & 1 & 1 \\
2 & 3 & 3 & 2 & 2 & 3 & 3 & 2 \\
3 & 2 & 2 & 3 & 3 & 2 & 2 & 3 \\
4 & 4 & 4 & 4 & 4 & 4 & 4 & 4 \\
\end{array}\right]}$ & No & $\mathbb{Z}_2$ &
$\displaystyle{\left[
\begin{array}{cccc|cccc}
1 & 1 & 1 & 1 & 1 & 1 & 1 & 1 \\
3 & 3 & 3 & 3 & 3 & 3 & 3 & 3 \\
2 & 2 & 2 & 2 & 2 & 2 & 2 & 2 \\
4 & 4 & 4 & 4 & 4 & 4 & 4 & 4 \\ \hline
1 & 4 & 4 & 1 & 1 & 4 & 4 & 1 \\
2 & 3 & 3 & 2 & 2 & 3 & 3 & 2 \\
3 & 2 & 2 & 3 & 3 & 2 & 2 & 3 \\
4 & 1 & 1 & 4 & 4 & 1 & 1 & 4 \\
\end{array}\right]}$ & No & $\mathbb{Z}_2\oplus \mathbb{Z}_2$ \\
& & & & & \\
$\displaystyle{\left[
\begin{array}{cccc|cccc}
2 & 2 & 1 & 1 & 2 & 2 & 1 & 1 \\
1 & 1 & 2 & 2 & 1 & 1 & 2 & 2 \\
4 & 4 & 4 & 4 & 4 & 4 & 4 & 4 \\
3 & 3 & 3 & 3 & 3 & 3 & 3 & 3 \\ \hline
2 & 2 & 1 & 1 & 2 & 2 & 1 & 1 \\
1 & 1 & 2 & 2 & 1 & 1 & 2 & 2 \\
3 & 3 & 4 & 4 & 3 & 3 & 4 & 4 \\
4 & 4 & 3 & 3 & 4 & 4 & 3 & 3 \\
\end{array}\right]}$ & No & $\mathbb{Z}_2\oplus \mathbb{Z}_2$ &
$\displaystyle{\left[
\begin{array}{cccc|cccc}
2 & 2 & 1 & 1 & 2 & 2 & 1 & 1 \\
1 & 1 & 2 & 2 & 1 & 1 & 2 & 2 \\
4 & 4 & 4 & 4 & 4 & 4 & 4 & 4 \\
3 & 3 & 3 & 3 & 3 & 3 & 3 & 3 \\ \hline
2 & 2 & 2 & 2 & 2 & 2 & 2 & 2 \\
1 & 1 & 1 & 1 & 1 & 1 & 1 & 1 \\
3 & 3 & 4 & 4 & 3 & 3 & 4 & 4 \\
4 & 4 & 3 & 3 & 4 & 4 & 3 & 3 \\
\end{array}\right]}$ & Yes & $\mathbb{Z}_2\oplus \mathbb{Z}_2$ \\
& & & & & \\
$\displaystyle{\left[
\begin{array}{cccc|cccc}
2 & 2 & 2 & 2 & 2 & 2 & 2 & 2 \\
1 & 1 & 1 & 1 & 1 & 1 & 1 & 1 \\
4 & 4 & 4 & 4 & 4 & 4 & 4 & 4 \\
3 & 3 & 3 & 3 & 3 & 3 & 3 & 3 \\ \hline
2 & 2 & 1 & 1 & 2 & 2 & 1 & 1 \\
1 & 1 & 2 & 2 & 1 & 1 & 2 & 2 \\
3 & 3 & 4 & 4 & 3 & 3 & 4 & 4 \\
4 & 4 & 3 & 3 & 4 & 4 & 3 & 3 \\
\end{array}\right]}$ & No & $\mathbb{Z}_4\oplus\mathbb{Z}_2$ &
$\displaystyle{\left[
\begin{array}{cccc|cccc}
2 & 2 & 2 & 2 & 2 & 2 & 2 & 2 \\
1 & 1 & 1 & 1 & 1 & 1 & 1 & 1 \\
4 & 4 & 4 & 4 & 4 & 4 & 4 & 4 \\
3 & 3 & 3 & 3 & 3 & 3 & 3 & 3 \\ \hline
2 & 2 & 2 & 2 & 2 & 2 & 2 & 2 \\
1 & 1 & 1 & 1 & 1 & 1 & 1 & 1 \\
3 & 3 & 4 & 4 & 3 & 3 & 4 & 4 \\
4 & 4 & 3 & 3 & 4 & 4 & 3 & 3 \\
\end{array}\right]}$ & No & $\mathbb{Z}_2\oplus \mathbb{Z}_2$ \\
& & & & & \\
\hline
\end{tabular}
\end{center}}
\caption{Non-quandle Biquandles of order 4 part 4.} 
\label{t6}
\end{table}

\begin{table} \footnotesize{
\begin{center}
\begin{tabular}{|ccc|ccc|}  \hline
  & & & & & \\
Biquandle & Self- & $\mathrm{Aut}(B)$ &
Biquandle & Self- & $\mathrm{Aut}(B)$\\
Matrix  & Flip?  & & Matrix & Flip? & \\
 & & & & & \\ \hline
 & & & & & \\
$\displaystyle{\left[
\begin{array}{cccc|cccc}
1 & 1 & 1 & 1 & 1 & 1 & 1 & 1 \\
3 & 2 & 2 & 3 & 3 & 2 & 2 & 3 \\
2 & 3 & 3 & 2 & 2 & 3 & 3 & 2 \\
4 & 4 & 4 & 4 & 4 & 4 & 4 & 4 \\ \hline
1 & 4 & 4 & 1 & 1 & 4 & 4 & 1 \\
2 & 2 & 2 & 2 & 2 & 2 & 2 & 2 \\
3 & 3 & 3 & 3 & 3 & 3 & 3 & 3 \\
4 & 1 & 1 & 4 & 4 & 1 & 1 & 4 \\
\end{array}\right]}$ & Yes & $\mathbb{Z}_2\oplus \mathbb{Z}_2$ &
$\displaystyle{\left[
\begin{array}{cccc|cccc}
1 & 1 & 2 & 2 & 1 & 1 & 2 & 2 \\
2 & 2 & 1 & 1 & 2 & 2 & 1 & 1 \\
4 & 4 & 4 & 4 & 4 & 4 & 4 & 4 \\
3 & 3 & 3 & 3 & 3 & 3 & 3 & 3 \\ \hline
1 & 1 & 1 & 1 & 1 & 1 & 1 & 1 \\
2 & 2 & 2 & 2 & 2 & 2 & 2 & 2 \\
3 & 3 & 4 & 4 & 3 & 3 & 4 & 4 \\
4 & 4 & 3 & 3 & 4 & 4 & 3 & 3 \\
\end{array}\right]}$ & No & $\mathbb{Z}_2\oplus \mathbb{Z}_2$ \\
& & & & & \\
$\displaystyle{\left[
\begin{array}{cccc|cccc}
1 & 1 & 2 & 2 & 1 & 1 & 2 & 2 \\
2 & 2 & 1 & 1 & 2 & 2 & 1 & 1 \\
4 & 4 & 4 & 4 & 4 & 4 & 4 & 4 \\
3 & 3 & 3 & 3 & 3 & 3 & 3 & 3 \\ \hline
1 & 1 & 2 & 2 & 1 & 1 & 2 & 2 \\
2 & 2 & 1 & 1 & 2 & 2 & 1 & 1 \\
3 & 3 & 4 & 4 & 3 & 3 & 4 & 4 \\
4 & 4 & 3 & 3 & 4 & 4 & 3 & 3 \\
\end{array}\right]}$ & No & $\mathbb{Z}_2\oplus \mathbb{Z}_2$ &
$\displaystyle{\left[
\begin{array}{cccc|cccc}
2 & 2 & 1 & 1 & 2 & 2 & 1 & 1 \\
1 & 1 & 2 & 2 & 1 & 1 & 2 & 2 \\
3 & 3 & 4 & 4 & 3 & 3 & 4 & 4 \\
4 & 4 & 3 & 3 & 4 & 4 & 3 & 3 \\ \hline
2 & 2 & 1 & 1 & 2 & 2 & 1 & 1 \\
1 & 1 & 2 & 2 & 1 & 1 & 2 & 2 \\
3 & 3 & 4 & 4 & 3 & 3 & 4 & 4 \\
4 & 4 & 3 & 3 & 4 & 4 & 3 & 3 \\
\end{array}\right]}$ & Yes & $\mathbb{Z}_4\oplus \mathbb{Z}_2$ \\
& & & & & \\
$\displaystyle{\left[
\begin{array}{cccc|cccc}
2 & 2 & 1 & 1 & 2 & 2 & 1 & 1 \\
1 & 1 & 2 & 2 & 1 & 1 & 2 & 2 \\
4 & 4 & 3 & 3 & 4 & 4 & 3 & 3 \\
3 & 3 & 4 & 4 & 3 & 3 & 4 & 4 \\ \hline
2 & 2 & 1 & 1 & 2 & 2 & 1 & 1 \\
1 & 1 & 2 & 2 & 1 & 1 & 2 & 2 \\
3 & 3 & 3 & 3 & 3 & 3 & 3 & 3 \\
4 & 4 & 4 & 4 & 4 & 4 & 4 & 4 \\
\end{array}\right]}$ & No & $\mathbb{Z}_4\oplus\mathbb{Z}_2$ &
$\displaystyle{\left[
\begin{array}{cccc|cccc}
2 & 2 & 2 & 2 & 2 & 2 & 2 & 2 \\
1 & 1 & 1 & 1 & 1 & 1 & 1 & 1 \\
3 & 3 & 4 & 4 & 3 & 3 & 4 & 4 \\
4 & 4 & 3 & 3 & 4 & 4 & 3 & 3 \\ \hline
2 & 2 & 2 & 2 & 2 & 2 & 2 & 2 \\
1 & 1 & 1 & 1 & 1 & 1 & 1 & 1 \\
3 & 3 & 4 & 4 & 3 & 3 & 4 & 4 \\
4 & 4 & 3 & 3 & 4 & 4 & 3 & 3 \\
\end{array}\right]}$ & Yes & $\mathbb{Z}_2\oplus \mathbb{Z}_2$ \\
& & & & & \\
$\displaystyle{\left[
\begin{array}{cccc|cccc}
2 & 2 & 2 & 2 & 2 & 2 & 2 & 2 \\
1 & 1 & 1 & 1 & 1 & 1 & 1 & 1 \\
4 & 4 & 3 & 3 & 4 & 4 & 3 & 3 \\
3 & 3 & 4 & 4 & 3 & 3 & 4 & 4 \\ \hline
2 & 2 & 2 & 2 & 2 & 2 & 2 & 2 \\
1 & 1 & 1 & 1 & 1 & 1 & 1 & 1 \\
3 & 3 & 3 & 3 & 3 & 3 & 3 & 3 \\
4 & 4 & 4 & 4 & 4 & 4 & 4 & 4 \\
\end{array}\right]}$ & No & $\mathbb{Z}_2\oplus\mathbb{Z}_2$ &
$\displaystyle{\left[
\begin{array}{cccc|cccc}
2 & 2 & 2 & 2 & 2 & 2 & 2 & 2 \\
1 & 1 & 1 & 1 & 1 & 1 & 1 & 1 \\
4 & 4 & 3 & 3 & 4 & 4 & 3 & 3 \\
3 & 3 & 4 & 4 & 3 & 3 & 4 & 4 \\ \hline
2 & 2 & 2 & 2 & 2 & 2 & 2 & 2 \\
1 & 1 & 1 & 1 & 1 & 1 & 1 & 1 \\
4 & 4 & 3 & 3 & 4 & 4 & 3 & 3 \\
3 & 3 & 4 & 4 & 3 & 3 & 4 & 4 \\
\end{array}\right]}$ & Yes & $\mathbb{Z}_2\oplus\mathbb{Z}_2$ \\
& & & & & \\
$\displaystyle{\left[
\begin{array}{cccc|cccc}
1 & 1 & 2 & 2 & 1 & 1 & 2 & 2 \\
2 & 2 & 1 & 1 & 2 & 2 & 1 & 1 \\
4 & 4 & 3 & 3 & 4 & 4 & 3 & 3 \\
3 & 3 & 4 & 4 & 3 & 3 & 4 & 4 \\ \hline
1 & 1 & 2 & 2 & 1 & 1 & 2 & 2 \\
2 & 2 & 1 & 1 & 2 & 2 & 1 & 1 \\
3 & 3 & 3 & 3 & 3 & 3 & 3 & 3 \\
4 & 4 & 4 & 4 & 4 & 4 & 4 & 4 \\
\end{array}\right]}$ & No & $\mathbb{Z}_2\oplus \mathbb{Z}_2$ &
$\displaystyle{\left[
\begin{array}{cccc|cccc}
1 & 1 & 1 & 1 & 1 & 1 & 1 & 1 \\
2 & 2 & 2 & 2 & 2 & 2 & 2 & 2 \\
4 & 3 & 3 & 3 & 4 & 3 & 3 & 3 \\
3 & 4 & 4 & 4 & 3 & 4 & 4 & 4 \\ \hline
1 & 1 & 1 & 1 & 1 & 1 & 1 & 1 \\
2 & 2 & 2 & 2 & 2 & 2 & 2 & 2 \\
4 & 3 & 3 & 3 & 4 & 3 & 3 & 3 \\
3 & 4 & 4 & 4 & 3 & 4 & 4 & 4 \\
\end{array}\right]}$ & Yes & $\mathbb{Z}_2$ \\
& & & & & \\
\hline
\end{tabular}
\end{center}}
\caption{Non-quandle Biquandles of order 4 part 5.} 
\label{t7}
\end{table}

\begin{table} \footnotesize{
\begin{center}
\begin{tabular}{|ccc|ccc|}  \hline
  & & & & & \\
Biquandle & Self- & $\mathrm{Aut}(B)$ &
Biquandle & Self- & $\mathrm{Aut}(B)$\\
Matrix  & Flip?  & & Matrix & Flip? & \\
 & & & & & \\ \hline
 & & & & & \\
$\displaystyle{\left[
\begin{array}{cccc|cccc}
1 & 1 & 1 & 1 & 1 & 1 & 1 & 1 \\
2 & 2 & 2 & 2 & 2 & 2 & 2 & 2 \\
4 & 4 & 3 & 3 & 4 & 4 & 3 & 3 \\
3 & 3 & 4 & 4 & 3 & 3 & 4 & 4 \\ \hline
1 & 1 & 1 & 1 & 1 & 1 & 1 & 1 \\
2 & 2 & 2 & 2 & 2 & 2 & 2 & 2 \\
4 & 4 & 3 & 3 & 4 & 4 & 3 & 3 \\
3 & 3 & 4 & 4 & 3 & 3 & 4 & 4 \\
\end{array}\right]}$ & Yes & $\mathbb{Z}_2\oplus \mathbb{Z}_2$ &
$\displaystyle{\left[
\begin{array}{cccc|cccc}
1 & 1 & 1 & 1 & 1 & 1 & 1 & 1 \\
2 & 2 & 2 & 2 & 2 & 2 & 2 & 2 \\
4 & 3 & 4 & 4 & 4 & 3 & 4 & 4 \\
3 & 4 & 3 & 3 & 3 & 4 & 3 & 3 \\ \hline
1 & 1 & 1 & 1 & 1 & 1 & 1 & 1 \\
2 & 2 & 2 & 2 & 2 & 2 & 2 & 2 \\
4 & 3 & 4 & 4 & 4 & 3 & 4 & 4 \\
3 & 4 & 3 & 3 & 3 & 4 & 3 & 3 \\
\end{array}\right]}$ & Yes & $\mathbb{Z}_2$ \\
& & & & & \\
$\displaystyle{\left[
\begin{array}{cccc|cccc}
1 & 1 & 1 & 1 & 1 & 1 & 1 & 1 \\
2 & 2 & 2 & 2 & 2 & 2 & 2 & 2 \\
4 & 4 & 4 & 4 & 4 & 4 & 4 & 4 \\
3 & 3 & 3 & 3 & 3 & 3 & 3 & 3 \\ \hline
1 & 1 & 1 & 1 & 1 & 1 & 1 & 1 \\
2 & 2 & 2 & 2 & 2 & 2 & 2 & 2 \\
4 & 4 & 4 & 4 & 4 & 4 & 4 & 4 \\
3 & 3 & 3 & 3 & 3 & 3 & 3 & 3 \\
\end{array}\right]}$ & Yes & $\mathbb{Z}_2\oplus\mathbb{Z}_2$ &
$\displaystyle{\left[
\begin{array}{cccc|cccc}
2 & 2 & 1 & 1 & 2 & 2 & 1 & 1 \\
1 & 1 & 2 & 2 & 1 & 1 & 2 & 2 \\
3 & 3 & 3 & 3 & 3 & 3 & 3 & 3 \\
4 & 4 & 4 & 4 & 4 & 4 & 4 & 4 \\ \hline
2 & 2 & 1 & 1 & 2 & 2 & 1 & 1 \\
1 & 1 & 2 & 2 & 1 & 1 & 2 & 2 \\
3 & 3 & 3 & 3 & 3 & 3 & 3 & 3 \\
4 & 4 & 4 & 4 & 4 & 4 & 4 & 4 \\
\end{array}\right]}$ & No & $\mathbb{Z}_2\oplus\mathbb{Z}_2$ \\
& & & & & \\
$\displaystyle{\left[
\begin{array}{cccc|cccc}
2 & 2 & 1 & 1 & 2 & 2 & 1 & 1 \\
1 & 1 & 2 & 2 & 1 & 1 & 2 & 2 \\
4 & 4 & 3 & 3 & 4 & 4 & 3 & 3 \\
3 & 3 & 4 & 4 & 3 & 3 & 4 & 4 \\ \hline
2 & 2 & 1 & 1 & 2 & 2 & 1 & 1 \\
1 & 1 & 2 & 2 & 1 & 1 & 2 & 2 \\
4 & 4 & 3 & 3 & 4 & 4 & 3 & 3 \\
3 & 3 & 4 & 4 & 3 & 3 & 4 & 4 \\
\end{array}\right]}$ & Yes & $\mathbb{Z}\oplus \mathbb{Z}_2$ &
$\displaystyle{\left[
\begin{array}{cccc|cccc}
2 & 2 & 2 & 2 & 2 & 2 & 2 & 2 \\
1 & 1 & 1 & 1 & 1 & 1 & 1 & 1 \\
4 & 4 & 4 & 4 & 4 & 4 & 4 & 4 \\
3 & 3 & 3 & 3 & 3 & 3 & 3 & 3 \\ \hline
2 & 2 & 2 & 2 & 2 & 2 & 2 & 2 \\
1 & 1 & 1 & 1 & 1 & 1 & 1 & 1 \\
4 & 4 & 4 & 4 & 4 & 4 & 4 & 4 \\
3 & 3 & 3 & 3 & 3 & 3 & 3 & 3 \\
\end{array}\right]}$ & Yes & $\mathbb{Z}_4\oplus\mathbb{Z}_2$ \\
& & & & & \\
$\displaystyle{\left[
\begin{array}{cccc|cccc}
1 & 1 & 2 & 2 & 1 & 1 & 2 & 2 \\
2 & 2 & 1 & 1 & 2 & 2 & 1 & 1 \\
4 & 4 & 3 & 3 & 4 & 4 & 3 & 3 \\
3 & 3 & 4 & 4 & 3 & 3 & 4 & 4 \\ \hline
1 & 1 & 2 & 2 & 1 & 1 & 2 & 2 \\
2 & 2 & 1 & 1 & 2 & 2 & 1 & 1 \\
4 & 4 & 3 & 3 & 4 & 4 & 3 & 3 \\
3 & 3 & 4 & 4 & 3 & 3 & 4 & 4 \\
\end{array}\right]}$ & Yes & $\mathbb{Z}_2\oplus\mathbb{Z}_4$ & & & \\
& & & & & \\
\hline
\end{tabular}
\end{center}}
\caption{Non-quandle Biquandles of order 4 part 6.} 
\label{t8}
\end{table}

\end{document}